%% file: source.tex
\documentclass[8pt,a4paper]{article}

\pdfoutput=1 

\usepackage[left=15mm,top=10mm,right=15mm,bottom=15mm]{geometry}
\usepackage{calc}
\usepackage{amsfonts}
\usepackage{float}
\usepackage{graphicx,epstopdf}
\usepackage{amsmath}
\usepackage[percent]{overpic}
\usepackage{multicol}
\usepackage{multirow}
\usepackage{algorithm}
\usepackage{dutchcal}
\usepackage[export]{adjustbox}
\usepackage{xargs} 
\usepackage[pdftex,dvipsnames]{xcolor}
\usepackage[colorinlistoftodos,prependcaption,textsize=scriptsize]{todonotes}
\usepackage{subfig}
\usepackage{graphicx,import}
\usepackage{caption}
\usepackage{gensymb}
\usepackage{soul,color}
\usepackage{authblk}

\newcommand{\mtx}[1]{\mathbf{#1}} 
\newcommand{\ten}[1]{\mathcal{#1}} 

\begin{document}

\title{Tensor B-Spline Numerical Methods for PDEs: a High-Performance Alternative to FEM}

\author[1]{Dmytro Shulga}
\author[2]{Oleksii Morozov}
\author[3]{Volker Roth}
\author[4]{Felix Friedrich}
\author[5]{Patrick Hunziker}

\affil[1]{University Hospital of Basel, Basel CH-4031, Switzerland (e-mail: shulgad@swissnano.org)}
\affil[2]{HighDim GmbH, Riehen CH-4125, Switzerland (e-mail: morozova@highdim.com)}
\affil[3]{University of Basel, Department of Mathematics and Computer Science, Basel CH-4051, Switzerland (e-mail: volker.roth@unibas.ch)}
\affil[4]{ETH Zurich, Department of Computer Science, Zurich CH-8092, Switzerland (email: felix.friedrich@inf.ethz.ch)}
\affil[5]{University Hospital of Basel, Clinic for Intensive Care, CLINAM, Basel CH-4031, Switzerland (e-mail: hunzikerp@swissnano.org)}

\date{}
\maketitle

\section*{Abstract}

Tensor B-spline methods are a high-performance alternative to solve partial differential equations (PDEs). This paper gives an overview on the principles of Tensor B-spline methodology, shows their use and analyzes their performance in application examples, and discusses its merits. Tensors preserve the dimensional structure of a discretized PDE, which makes it possible to develop highly efficient computational solvers. B-splines provide high-quality approximations, lead to a sparse structure of the system operator represented by shift-invariant separable kernels in the domain, and are mesh-free by construction. Further, high-order bases can easily be constructed from B-splines. 
In order to demonstrate the advantageous numerical performance of tensor B-spline methods, we studied the solution of a large-scale heat-equation problem (consisting of roughly 0.8 billion nodes!) on a heterogeneous workstation consisting of multi-core CPU and GPUs. Our experimental results nicely confirm the excellent numerical approximation properties of tensor B-splines, and their unique combination of high computational efficiency and low memory consumption, thereby showing huge improvements over standard finite-element methods (FEM). 

\bigskip
\textbf{Keywords:} partial differential equation, Ritz-Galerkin formulation, tensor algebra, B-spline, finite-element method, shift invariance, filter, kernel, structure, parallel processing, sparse matrix-vector multiplication, non-uniform memory access, multi-core processor.

\section{Introduction}

Finding a numerical solution of a PDE is a computationally intensive task in many scientific and engineering applications. Even ``routine'' modeling problems such as weather forecasting can be extremely challenging in practice, and the  development of fast and accurate PDE solvers is an important field of research.

Finite-element numerical methods (FEM) for PDEs rely on mesh-based domain discretization and employ polynomial basis functions. FEM has been studied extensively for decades and became the de-facto instrument for solving PDEs on arbitrarily-shaped domains. Meshing automation, however, defines a key practical problem in FEM-based approaches, which heavily depends on the specific properties of the domain considered \cite{bathe2019}. While pure mesh-less methods require much less effort to discretize the domain, finding efficient numerical integration schemes is challenging for such methods \cite{De2000}. Discontinuous Galerkin (DG) FEM methods combine physical accuracy and flexibility of mesh generation via weakly enforced continuity of discontinuous elements, however at the price of relatively high computational costs \cite{Kaufmann2009}.

Contrary to FEM-based approaches characterized by polynomial bases and the use of matrix algebra,  the class of tensor B-spline methods discussed in this work uses spline bases and tensors. Splines have been used for computer graphics and computer-aided design for a long time, and they also have been applied to signal 
and image processing and reconstruction.  A specific B-spline framework for signal and image processing was proposed by Unser  \cite{Unser1999, Unser1993part1, Unser1993part2}. The use of splines for solving PDEs has been studied by H{\"o}llig \cite{Hollig2003, Hollig2005, Hollig2015}, who showed how B-splines can be used in the context of finite-element methods. The classical approach to solve the resulting discrete system of equations is based on matrix algebra using sparse matrix formats. However, this approach suffers from high computational costs which limits its applicability to problems of small or moderate size. To some extent, these computational problems can be overcome in a very elegant way by using tensors \cite{Cichocki2015, Sidiropoulos2017}. Compared to matrices, such tensors allow us to represent multidimensional structures in a more compact and natural way. This computational tensor algebra approach  \cite{Morozov2010}  makes it possible to develop highly efficient numerical algorithms, and several benefits of combining tensors and B-splines have been shown in the context of multidimensional signal reconstruction \cite{Morozov2011} and for solving diffusion PDEs in Optical Diffusion Tomography \cite{Shulga2017, Shulga2018}.

The Tensor B-spline method has many appealing properties of a ``generic'' numerical PDE solver: 1) it provides us with accurate and flexible discretizations of coefficients, sources and solutions, 2) it allows for efficient integration strategies, 3) it makes it relatively simple to develop fast and memory-efficient algorithms, 4) the mathematical elegance of computational tensor algebra \cite{Morozov2010} leads to natural and transparent models. The key advantages over FEM are: 1) no mesh is needed, 2) a high-degree B-spline is more efficient than a high-degree FEM polynomial, and 3) a high-degree Tensor B-spline solver is more computationally and memory efficient than classical FEM solvers. Therefore, splines offer accurate solutions at significantly reduced computational costs.

\begin{figure}[H]
\centering     
\includegraphics[scale=0.55]{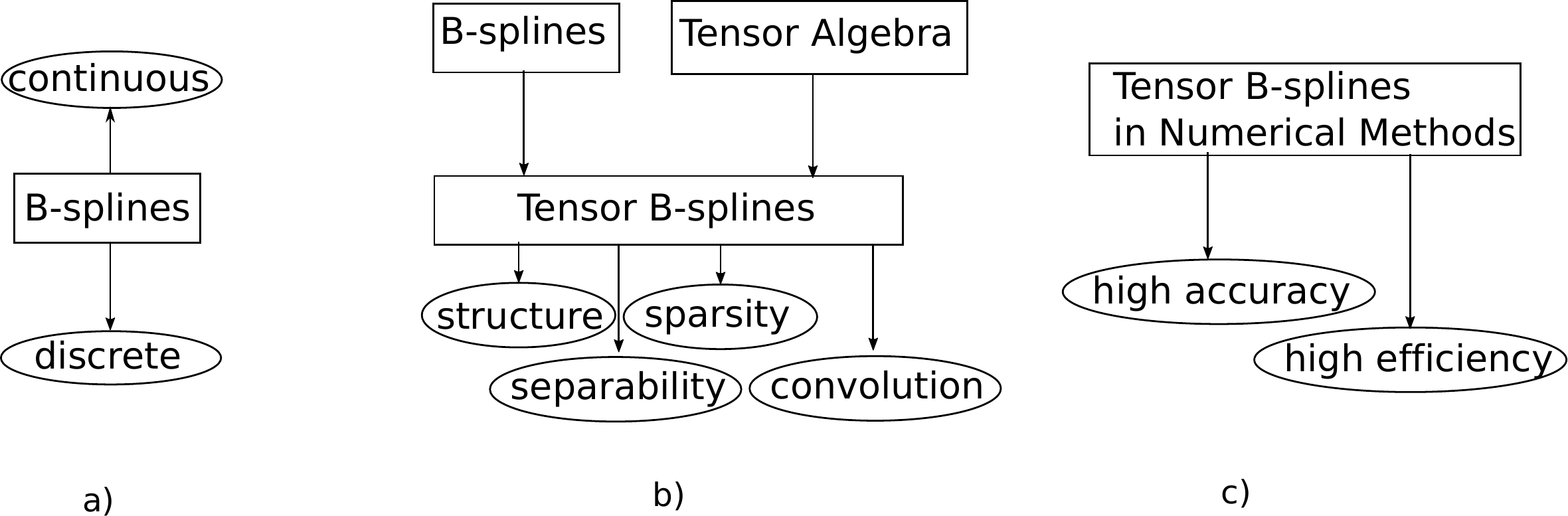}
\caption{Diagram of Tensor B-spline PDE solver main features} \label{figTensBspline}
\end{figure}

Fig.~\ref{figTensBspline} highlights the main properties that make Tensor B-splines promising candidates for numerical PDE methods. Splines naturally link continuous 
and discrete domains \cite{Unser1999} (Fig.~\ref{figTensBspline} (a)) while providing excellent approximations of coefficients, sources and solutions. The combination of B-splines and 
tensor algebra preserves the intrinsic structure of the problem, and enforces both sparsity and separability in a very natural way (Fig.~\ref{figTensBspline} (b)). At the same time, it makes it easy to design highly efficient parallel and matrix-free algorithms. As a result, highly accurate and efficient solutions can be obtained (Fig.~\ref{figTensBspline} (c)).

The remainder of this paper is structured as follows.
Section \ref{secTensBsplineMethod} introduces the mathematical concepts behind Tensor B-spline numerical methods for PDEs.  Elliptic PDEs, B-spline spaces,  and fundamental properties of  B-splines are discussed, with specific emphasis on the use of computational tensor algebra to the B-spline Ritz-Galerkin formulation. Further, efficient computational strategies are discussed. Section \ref{secMethodPerf} focuses on assessing the performance of the method proposed based on several real-world examples. Finally, Sections \ref{secDiscussion} and \ref{secConclusion} present an in-depth discussion and conclusion.

\section{Tensor B-spline Method for PDEs} \label{secTensBsplineMethod}

\subsection{Notation}

\begin{figure}[H]
\centering
\def\svgscale{0.7}
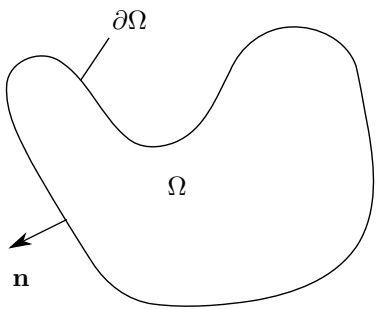
\caption{A domain $\Omega$, the domain boundary $\partial \Omega$, and a normal to the domain boundary $\mathbf{n}$}  \label{figDomain}
\end{figure}

Vectors are denoted by lower-case bold symbols, such as $\mathbf{b} \in \mathbb{R}^{d}$. Particularly we denote $\mathbf{x} \in \mathbb{R}^{d},~ d=1,2,3$ 
to be an element of real coordinate space. An element of a vector $\mathbf{c}$ is denoted as $\mathbf{c}_k$. Matrices are denoted by upper-case bold symbols, 
such as $\mathbf{A} \in \mathbb{R}^{d_1 \times d_2}$. If not specified differently, scalar fields are denoted as $ D(\mathbf{x}) \in \mathbb{R}$.
Tensors are denoted by calligraphic symbols \cite{Morozov2010} with shorthand notation for multiple tensor indices, such as $\mathcal{D}^{\mathbf{m}}$, 
$\mathcal{w}_{\mathbf{lm}} \in \mathbb{R}$, where $\mathbf{l}, ~ \mathbf{m} \in \mathbb{Z}^{d}, ~d=1,2,3$. We define $\Omega$ to be a domain with boundary 
$\partial \Omega$ and normal to a domain boundary $\mathbf{n}$ (see Fig.~\ref{figDomain}).

\subsection{Elliptic PDEs}

A generic partial differential equation (PDE)
\begin{eqnarray}
F(\mathbf{x},\varphi(\mathbf{x}),  \varphi'_{x_{1}}, ..., \varphi'_{x_{n}},   \varphi''_{x_{1}}, ..., \varphi''_{x_{n}}, ...) = 0, \mathbf{x} \in \Omega 
\end{eqnarray}
links an unknown multivariate function $\varphi(\mathbf{x})$ and its partial derivatives in the domain $\Omega$. In order to obtain a solution of the PDE,  
boundary conditions (BC) have to be provided to define the function behavior on the domain boundary $\partial \Omega$.

In most cases it is either impossible or intractable to obtain an analytic solution of a PDE coupled with BC on domains $\Omega$ of arbitrary shape.  
Therefore the solution is obtained numerically. The unknown function is approximated by an expansion 
\begin{equation}
\varphi(\mathbf{x}) \approx \hat{\varphi}(\mathbf{x}) = \sum_{k \in \mathbb{Z}} c_{k} \eta(\mathbf{x}), ~\mathbf{x} \in \Omega.  \label{expansion}
\end{equation} 
From Equation \eqref{expansion} it follows that for given basis functions $\eta(\mathbf{x})$, the function $\hat{\varphi}(\mathbf{x})$ itself is fully 
described by the expansion coefficients $c_{k}$. In practice, it is highly important to choose appropriate basis functions $\eta(\mathbf{x})$, where ``appropriate'' typically refers to good approximation properties and linear independence. The error between the function $\varphi(\mathbf{x})$ and its approximation 
$\hat{\varphi}(\mathbf{x})$  defines the residual, which is then minimized by specific numerical procedures. For instance, the method of weighted residuals requires  
\begin{equation}
\int_{\Omega} \psi_{l}(\mathbf{x}) (\varphi(\mathbf{x}) - \hat{\varphi}(\mathbf{x})) d\mathbf{x} = 0, ~ l \in \mathbb{Z}.  \label{mwr}
\end{equation}
If the weight functions are chosen to be equal to basis functions $\psi_{l}(\mathbf{x}) = \eta_{l}(\mathbf{x})$, the Ritz-Galerkin formulation is obtained
(which is equivalent to the least squares method). Substitution of \eqref{expansion} into \eqref{mwr} leads to 
$\sum_{k}c_{k} \int_{\Omega} \psi_{l} \eta_{k}(\mathbf{x})d\mathbf{x} = \int_{\Omega} \psi_{l} \varphi(\mathbf{x}) d\mathbf{x},~ l \in \mathbb{Z}$ and can be written
in terms of the variational formulation
\begin{equation}
a( \varphi(\mathbf{x}), \psi(\mathbf{x}) ) = l(\psi(\mathbf{x})), 
\end{equation}
where $a(\cdot,\cdot)$ is an elliptic bilinear form and $l(\cdot)$  is a bounded linear functional on a Hilbert space $\mathbb{H}$. Basis (trial) functions
$~\varphi(\mathbf{x})$ and weight (test) functions $\psi(\mathbf{x})$  belong to Sobolev space: 
$\varphi(\mathbf{x})  \in  \mathbb{H}^{1}(\Omega)~ ,\psi(\mathbf{x})  \in  \mathbb{H}^{1}(\Omega),~\mathbb{H}^{1}(\Omega) = \lbrace f(x):  ||f(x)||_{\mathbb{L}^{2}(\Omega)} < \infty, ~ ||f(x)||_{\mathbb{H}^{1}(\Omega)} < \infty \rbrace$ 
(see Appendix I for the exact definitions of the norms).

We now consider elliptic PDEs of the form
\begin{eqnarray}\label{elliptic_pde}
- \nabla \cdot (D(\mathbf{x})\nabla \varphi(\mathbf{x})) + \mu_{a}(\mathbf{x}) \varphi(\mathbf{x}) = q(\mathbf{x}), ~\mathbf{x} \in ~ \Omega. \label{diffpde} 
\end{eqnarray}
If we assume that \eqref{elliptic_pde} describes some diffusion process, then $\varphi(\mathbf{x})$ is the density of the diffusing material, $D(\mathbf{x})$ is 
the diffusion coefficient,  $\mu_a(\mathbf{x})$ is  the absorption coefficient and $q(\mathbf{x})$ is the source density. As special cases we obtain the Poisson equation 
$ \nabla^{2} \varphi(\mathbf{x}) = q(\mathbf{x}),~\mathbf{x} \in \Omega$ and the Laplace equation $\nabla^{2} \varphi(\mathbf{x}) = 0, ~\mathbf{x}\in \Omega$.

An elliptic PDE can be coupled with boundary conditions of different types. For $\mathbf{x} \in \partial \Omega$: 1) Dirichlet BC $\varphi(\mathbf{x}) = g(\mathbf{x})$ (non-homogeneous), $\varphi(\mathbf{x}) = 0$ (homogeneous);
2) Neumann BC $\nabla \varphi(\mathbf{x}) \cdot \mathbf{n} = g(\mathbf{x})$, $\mathbf{n}$ is the outward normal to the domain boundary $\partial \Omega$;
3) Robin BC $\alpha(\mathbf{x}) (\nabla \varphi(\mathbf{x}) \cdot \mathbf{n}) + \beta(\mathbf{x})\varphi(\mathbf{x}) = g(\mathbf{x})$; 
4) Cauchy BC $\varphi(\mathbf{x}) = a(\mathbf{x})$, $\nabla \varphi(\mathbf{x}) \cdot \mathbf{n} = b(\mathbf{x})$;
5) Mixed BC requires different boundary conditions to be satisfied on disjoint parts of the boundary of the domain where the condition is stated. 

In the case of the Robin BC of the form $2D(\mathbf{x}) (\nabla \varphi(\mathbf{x}) \cdot \mathbf{n}) + \varphi(\mathbf{x}) = 0, \mathbf{x} \in \partial \Omega$, 
after in integration \eqref{diffpde} by parts we have
\begin{eqnarray}
a(\varphi(\mathbf{x}), \psi(\mathbf{x})) = \int_{\Omega} D(\mathbf{x}) \nabla \varphi(\mathbf{x}) \cdot \nabla \psi(\mathbf{x}) + \mu_{a}(\mathbf{x}) \varphi(\mathbf{x}) \psi(\mathbf{x}) d \mathbf{x} + \frac{1}{2} \int_{\partial \Omega} \varphi(\mathbf{x}) \psi(\mathbf{x}) d s,  \label{a} \\
l(\psi(\mathbf{x})) = \int_{\Omega} q(\mathbf{x}) \psi(\mathbf{x}) d\mathbf{x}. \label{l}
\end{eqnarray}

The important idea underlying finite element methods (FEMs) is the specific choice of basis functions $\eta(\mathbf{x})$ with limited (local) support. 
For example, FEM uses Lagrange polynomials that are non-zero only within an element, and equal to zero outside. In the next section we show how the unknown function
$\varphi(\mathbf{x})$ can be represented in B-spline basis. 

\subsection{B-spline Spaces}

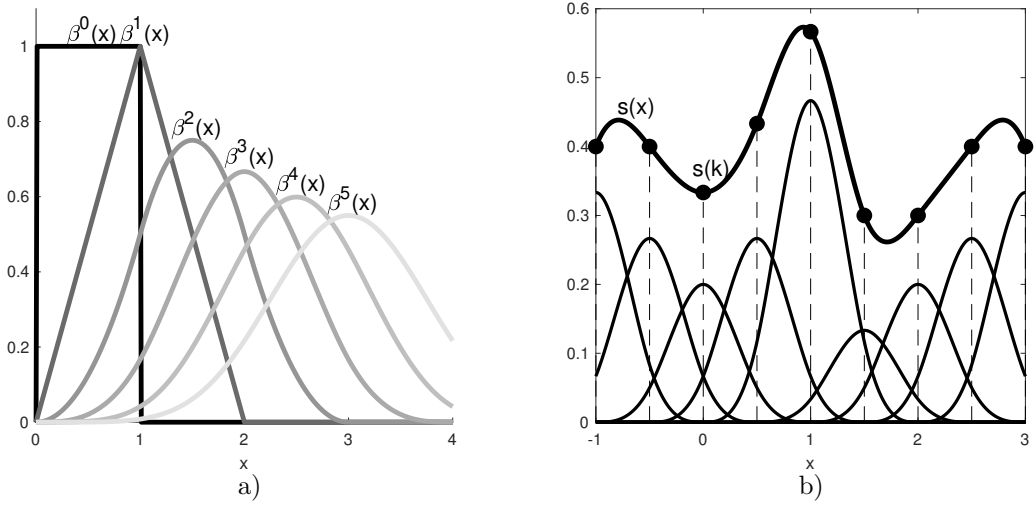
\begin{figure}[H]
\centering
    \includegraphics[scale=0.6]{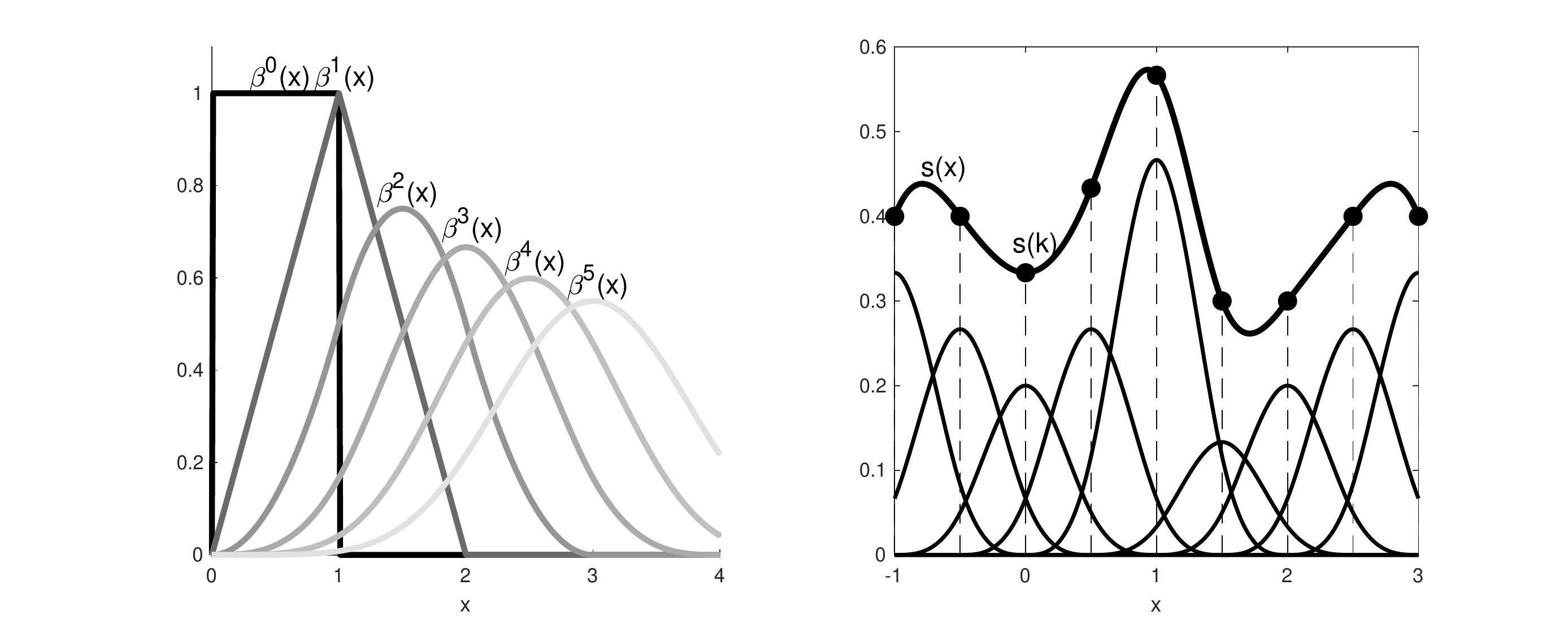}
    \begin{picture}(400,15)
\put(100,10){a)}
\put(310,10){b)}
\end{picture}
    \caption{(a) Univariate B-spline functions $\beta^{n}(x)$ of degree $n=0...5$, (b) An example of interpolation of samples $s(k)$ using univariate cubic
    B-spline functions $\beta^{3}(x)$, $s(x) = \sum_{k} c(k) \beta^{3}(x-k).$} \label{figBsplines}
\end{figure}

Before we apply B-spline basis functions to \eqref{a} and \eqref{l}, we would like to give a short overview of their properties. Since their 
introduction in the late 60's \cite{schoenberg1969, deboor1972}, B-splines have found many applications in computer graphics, computer-aided design, medical imaging \cite{Unser1999, Unser1993part1, Unser1993part2}, 
PDEs \cite{Hollig2003, Hollig2005, Hollig2015}, etc. A graphical representation of univariate B-spline functions is shown in Fig.~\ref{figBsplines} (a), B-splines of higher 
dimensions (multivariate B-splines) are obtained via tensor products. An example of a one-dimensional interpolation using cubic B-splines is shown in Fig.~\ref{figBsplines} (b). 

We consider B-splines on uniform grids. Table~\ref{tblBsplines} provides a summary of important properties of B-splines. For a detailed theory of B-splines we refer the reader to \cite{Hollig2003, Unser1999}.

\begin{table}[H]
\centering
\def\arraystretch{1.6}
\begin{tabular}{ l  c  c }
\hline
sign and support &  positive in its local support $(-\frac{n+1}{2}, \frac{n+1}{2})$, zero outside  \\
smoothness &   $(n-1)$-times continuously differentiable \\
structure &  piecewise polynomial with smoothly connected pieces, symmetrical,  monotone\\
derivative &  $\frac{d \beta^{n}(x)}{dx} = \beta^{n-1}(x+\frac{1}{2}) - \beta^{n-1}(x-\frac{1}{2}) $ \\
integration &  $\int_{-\infty}^{x} \beta^{n}(x)dx = \sum_{k=0}^{+\infty} \beta^{n+1} (x-\frac{1}{2}-k)$\\
convolution & $\beta^{m+n+1}(x) = \int_{supp} \beta^{m}(x-y)\beta^{n}(y)dy$,~  $\beta^{n}(x) = \underbrace{\beta^{0}(x) \ast \beta^{0}(x) \ast ... \ast \beta^{0}(x)}_{(n+1)~\text{times}} $ \\
scalar product &  $\beta^{m+n+1}(k-l) = \int_{supp} \beta^{m}(x-k) \beta^{n}(x-l) dx$\\ 
interpolation &  $s(x) = \sum_{k \in \mathbb{Z}} c_k \beta^{n}(x-k)$, $c_k =(b^{n})^{-1} \ast s_k$,
cubic spline has minimum curvature property \\
cardinal representation &  $s(x) = \sum_{k \in \mathbb{Z}} s_k \eta^{n}(x-k)$, $\eta^{n}(x) = \sum_{k \in \mathbb{Z}} (b^{n}_{k})^{-1} \beta^{n}(x-k)$ \\
\hline
\end{tabular}
\caption{The summary of  properties  of a univariate B-spline $\beta^{n}(x)$.} \label{tblBsplines}
\end{table}

\begin{figure}[H]
\scriptsize
\centering
\def\svgscale{0.6}
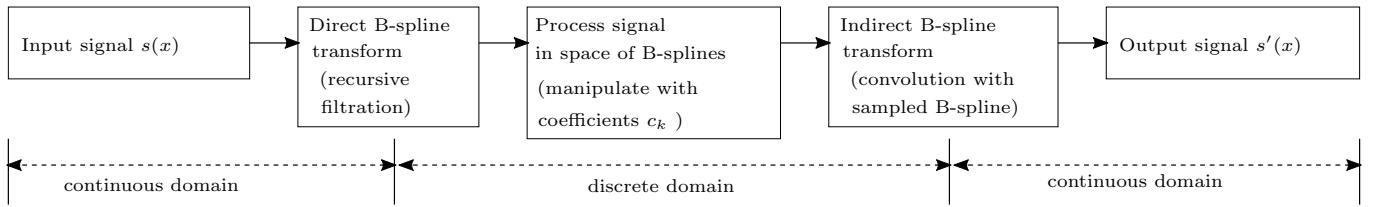
\caption{Signal processing in the B-spline framework}  \label{figBspnSignProc}
\normalsize
\end{figure}

A generic framework for signal processing with B-splines was extensively studied  by Unser  \cite{Unser1993part1, Unser1993part2}. The basic steps are shown in Fig.~\ref{figBspnSignProc}:
1) the continuous input signal $s(x)$ is transformed into the B-spline space via so-called direct B-spline transform (implemented as recursive filtration); 2) the signal processing is performed in the discrete domain of B-spline coefficients $c_{k}$; 3) the result is transferred back to the continuous signal domain via so-called indirect B-spline transform (implemented as convolution with sampled B-spline). 

The B-spline framework utilizes efficient digital filtering techniques to perform the interpolation, approximation, differentiation, etc.\ of the multidimensional input signal. This approach to digital filtering has been shown to be more efficient than matrix-based algorithms \cite{Unser1993part2}.

\begin{figure}[H]
\centering     
\includegraphics[scale=0.8]{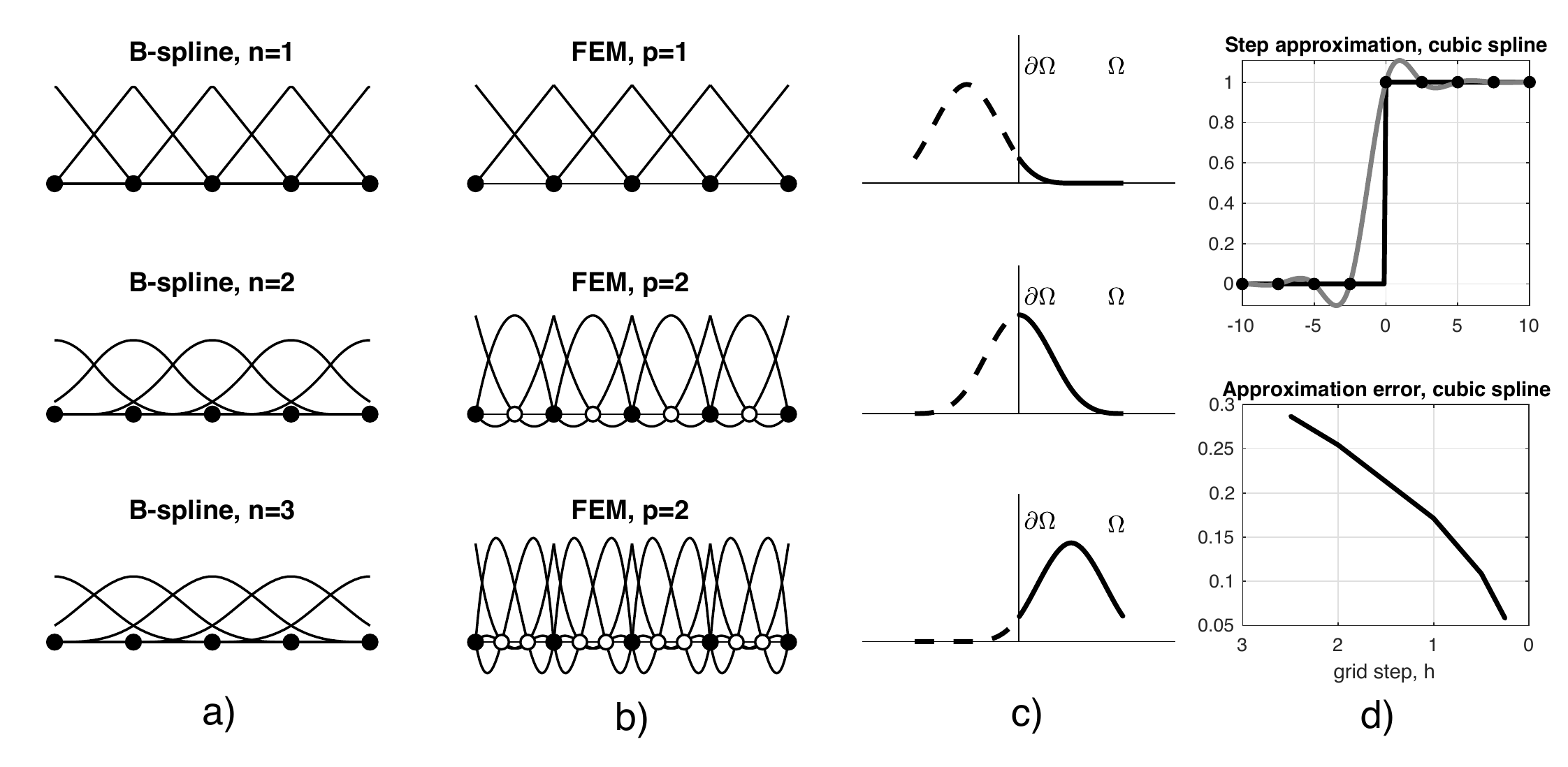}
\caption{Univariate B-spline (a) and FEM (b) bases,  an intersection of a B-spine function and a domain boundary (c)} \label{figFEMSplineBasis}
\end{figure}

We depict the univariate B-splines in comparison to the FEM polynomial functions in Fig.~\ref{figFEMSplineBasis} (a, b). 

The one-dimensional example clearly shows the main differences between the bases. While linear FEM and B-splines coincide, an increasing basis order introduces additional nodes in FEM (white circles) and spreads the support of B-splines. B-splines do not conform to the domain boundary (c): for $n>1$, B-splines beyond the domain also contribute to the solution. Fig.~\ref{figFEMSplineBasis} (d) shows the approximation of the step function with cubic B-spline and the approximation error. The reader might notice that, when the step function is sampled, the B-splines coinciding with grid nodes will result in an exact representation of the step function.

A multivariate B-spline of degree $n$ with grid step $\mathbf{h} \in \mathbb{R}^d$ is defined as outer (tensor) product of $d$ univariate B-splines:
\begin{equation}
\beta^{n}_{\mathbf{k},\mathbf{h}} (\mathbf{x}) = \beta^{n} (x_{1}/h_{1} - k_{1}) \cdots  \beta^{n} (x_{d}/h_{d} - k_{d}) ,~ \mathbf{k} \in \mathbb{Z}^{d}
\end{equation}

By analogy to \eqref{expansion} we expand unknown function (solution of a PDE) by B-spline basis functions

\begin{align}
\hat{\varphi}(\mathbf{x}) &= \sum_{k_1} \sum_{k_2} \cdots \sum_{k_d} \mathcal{c}_{k_1 k_2 \cdots k_d}  \beta^{n} (x_{1}/h_{1} - k_{1}) \cdots  \beta^{n} (x_{d}/h_{d} - k_{d}) \nonumber  \\ 
										&= \sum_{\mathbf{k} \in \mathbb{Z}^{d}} \mathcal{c}_{\mathbf{k}} \beta^{n}_{\mathbf{k},\mathbf{h}} (\mathbf{x}). 
\end{align}
The same operation is applied to the known coefficients and source functions:
\begin{equation}
D(\mathbf{x}) = \sum_{\mathbf{m} \in \mathbb{Z}^{d}} \mathcal{d}_{\mathbf{m}} \beta^{n}_{\mathbf{m},\mathbf{h}} (\mathbf{x}), ~
\mu_{a}(\mathbf{x}) = \sum_{\mathbf{m} \in \mathbb{Z}^{d}} \mathcal{d}_{\mathbf{m}} \beta^{n}_{\mathbf{m},\mathbf{h}} (\mathbf{x}), ~
q(\mathbf{x}) = \sum_{\mathbf{j} \in \mathbb{Z}^{d}} \mathcal{q}_{\mathbf{j}} \beta^{n}_{\mathbf{j},\mathbf{h}} (\mathbf{x}).
\end{equation}
The expansion coefficients $\mathcal{d}_{\mathbf{m}},~\mathcal{m}_{\mathbf{m}},~\mathcal{q}_{\mathbf{j}}$ are obtained via direct B-spline transform that exploits either 
interpolation of approximation ($\mathbb{L}^{2}$ projection); after solving  for $\mathcal{c}_{\mathbf{k}}$ they have to be transformed to the signal 
space using indirect B-spline transform (see Fig.~\ref{figBspnSignProc}). 

Finally, we have a discrete Ritz-Galerkin formulation of \eqref{a}, \eqref{l}:
\begin{eqnarray} \label{RitzGalerkin}
\sum_{\mtx{k}} \ten{c}_{\mtx{k}} \sum_{\mtx{m}} \ten{d}_{\mtx{m}} \int_{\Omega}  
\left( \nabla \beta^{n_{b}}_{\mtx{k},\mtx{h}}(\mtx{x})  \cdot \nabla \beta^{n_{b}}_{\mtx{l},\mtx{h}}(\mtx{x}) \right)  \beta^{n_{p}}_{\mtx{m},\mtx{h}}(\mtx{x}) d \mtx{x}  \nonumber 
+\sum_{\mtx{k}} \ten{c}_{\mtx{k}} \sum_{\mtx{m}} \ten{m}_{\mtx{m}} \int_{\Omega}  
\beta^{n_{b}}_{\mtx{k},\mtx{h}}(\mtx{x})  \beta^{n_{b}}_{\mtx{l},\mtx{h}}(\mtx{x})  \beta^{n_{p}}_{\mtx{m},\mtx{h}}(\mtx{x}) d \mtx{x}  \nonumber \\ 
+\frac{1}{2\gamma} \sum_{\mtx{k}} \ten{c}_{\mtx{k}} \int_{\partial \Omega}\! \beta^{n_{b}}_{\mtx{k},\mtx{h}}(\mtx{x})  \beta^{n_{b}}_{\mtx{l},\mtx{h}}(\mtx{x}) ds  
=\sum_{\mtx{j}} \ten{q}_{\mtx{j}} \int_{\Omega} \beta^{n_{s}}_{\mtx{j},\mtx{h}}(\mtx{x})  \beta^{n_{b}}_{\mtx{l},\mtx{h}}(\mtx{x}) d \mtx{x}.
\end{eqnarray}

The rest of the paper is dedicated to  efficient solution strategies for Equation \eqref{RitzGalerkin} or similar formulations involving the use of tensor algebra and efficient filter-like numerical algorithms, allowing  efficient parallel implementations for multi-core processors and GPUs. 

\subsection{Tensor Structure}

The discrete Ritz-Galerkin formulation \eqref{RitzGalerkin} is inherently multidimensional. Indeed, the formulation has 
more dimensions then the initial PDE problem. To deal with this situation, two main ideas are commonly used: 
1) the basis functions are chosen with a small support in order to make the discretization sparse and therefore the problem computationally feasible 
(crucial idea of FEM);
2) the multidimensional formulation is folded into sparse matrices and vectors to fit well-established routines of Matrix 
Algebra \cite{Hollig2003, Hollig2005, Zienkiewicz1977, Smith2013}, and afterwards the solution is rearranged into the original dimensions of the problem.

However, the standard approach of matricizing the multidimensional formulation has its limitations. While it flattens and merges the different dimensions, the underlying structure (containing important information for efficient computations) appears to be hidden. Given this flattened representation, there is only limited room for optimization, mainly dealing with values and indices of the block-diagonal sparse matrix format. This format, however, typically has little in common with the initial problem structure. 
Moreover, the structure of the (sparse) matrix needs to be represented, adding overhead to the implementation and rendering a software framework less generic.

\begin{figure}[H]
\center
\normalsize
\def\svgscale{1.0}
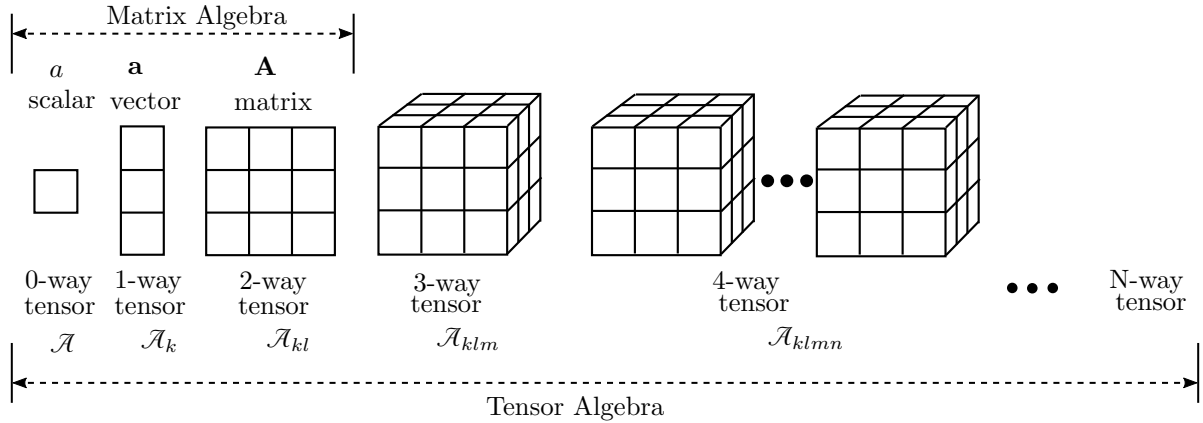
\caption{Tensors as generalizations of scalars, vectors and matrices.} \label{figTensor}
\normalsize
\end{figure}

Tensors, or multi-way arrays, are natural objects to be used instead of matrices in such cases. Tensors are replacing matrices more and more in many problems that were originally described in terms of matrices \cite{Kolda2009, Cichocki2015, Sidiropoulos2017}. From a technical perspective,  tensors are straight-forward  generalizations of vectors and matrices, cf.~Fig.~\ref{figTensor}. 

Tensors preserve the dimensional structure and data coherence. While being slightly more complicated objects than matrices, 
tensors frequently allow for more elegant solutions in a simpler way. The natural instrument for solving tensor structured problems is Computational Tensor Algebra \cite{Morozov2010}, targeted to highlight
possible optimizations for computations with multidimensional data.  The summary of tensor operations in comparison to operations with matrices and vectors is presented in Table~\ref{tblOp}. As one can see, tensor algebra operations are fully defined by tensor indices and their positions rather than by
special symbols like, for example, the Kronecker product in matrix algebra. Such a unification has the potential to allow for automatic simplification of tensor expression and optimal derivation of computational algorithms. For details we refer the reader to \cite{Morozov2010}.

\begin{table}[H]
\centering
\def\arraystretch{1}
\begin{tabular}{ l  c  c }
\hline
Operation                     & Tensor Algebra                                                     & Matrix Algebra  \\
\hline
inner product                 & $a = \mathcal{u}_{n}$ $\mathcal{v}^{n}$                            & $a = \mathbf{u}^{T} \mathbf{v}$\\
outer product                 & $\mathcal{A}^{n_{1}n_{2}} = \mathcal{u}^{n_1}$ $\mathcal{v}^{n_2}$ & $A = \mathbf{u} \mathbf{v}^{T}$\\
matrix-vector product         & $\mathcal{v}^{m} = \mathcal{B}^{m}_{n}\mathcal{u}^{n}$             & $\mathbf{v} = \mathbf{B} \mathbf{u}$ \\
matrix-matrix product         & $\mathcal{D}^{m}_{l} = \mathcal{B}^{m}_{n}\mathcal{C}^{n}_{l}$     & $\mathbf{D} = \mathbf{B} \mathbf{C}$\\
element-by-element product    &$\mathcal{w}_{m} = \mathcal{u}^{m}\mathcal{v}^{m}$                  & $\mathbf{w} = \mathbf{u} \circ \mathbf{v}$\\
tensor product                & $\mathcal{F}^{mn}_{kl} = \mathcal{B}^{mn}\mathcal{E}_{kl}$         & $\mathbf{F} = \mathbf{B} \otimes \mathbf{E}$\\								
\hline
\end{tabular}
\caption{The summary of important mathematical operations in Matrix and Tensor Algebra, where $a \in \mathbb{R}$ is a scalar, $\mathbf{u},~\mathbf{v},~\mathbf{w} \in \mathbb{R}^{n}$ are vectors, $\mathbf{A}\in \mathbb{R}^{n \times n},~ \mathbf{B}  \in \mathbb{R}^{m \times n},~  \mathbf{C}  \in \mathbb{R}^{n \times l},~ \mathbf{D} \in \mathbb{R}^{m \times \l},~ \mathbf{E} \in \mathbb{R}^{k \times \l},~ \mathbf{F} \in \mathbb{R}^{mk \times nl}$ are matrices, $(\cdot)^{T}$ - transposition operation, $\circ$ - Hadamard product, $\otimes$ - Kroneker product, $\mathcal{A}\in \mathbb{R}^{n \times n},~ \mathcal{B}  \in \mathbb{R}^{m \times n},~  \mathcal{C}  \in \mathbb{R}^{n \times l},~ \mathcal{D} \in \mathbb{R}^{m \times \l},~ \mathcal{E} \in \mathbb{R}^{k \times \l},~ \mathcal{F} \in \mathbb{R}^{m \times n \times k \times l}$ are tensors.} \label{tblOp}
\end{table}

We will use the described tensor notation further in this paper. After applying the symbols of tensor algebra to the Ritz-Galerkin formulation \eqref{RitzGalerkin} we get
\begin{eqnarray}  
\mathcal{c}_{\mathbf{k}} \mathcal{d}_{\mathbf{m}} \mathcal{w}^{\mathbf{klm}} + 
\mathcal{c}_{\mathbf{k}} \mathcal{m}_{\mathbf{m}} \mathcal{f}^{\mathbf{klm}} + 
\frac{1}{2} \mathcal{c}_{\mathbf{k}} \mathcal{h}^{\mathbf{kl}} = \mathcal{q}_{\mathbf{j}} \mathcal{r}^{\mathbf{jl}} 
~~\Leftrightarrow ~~
F(\mathcal{c}) = \mathcal{t},\label{tensForm}
\end{eqnarray}
where tensors $\mathcal{w}^{\mathbf{klm}},~\mathcal{f}^{\mathbf{klm}},~\mathcal{h}^{\mathbf{kl}},~\mathcal{r}^{\mathbf{ji}}$ correspond to integrals in \eqref{RitzGalerkin}. In expression \eqref{tensForm}, the multidimensional integrals are encapsulated, and arithmetic operations can be considered in terms of tensor algebra indices as shown in Table~\ref{tblOp}. First, the expression suggests that computations can be done either via an inner or an outer product. Second, it can be observed that different algorithms are defined where one of the indices $\mathbf{k}$, $\mathbf{l}$ or $\mathbf{m}$ is used in the algorithm's outermost loop. 

When  the system of equations \eqref{tensForm} has a large number of unknowns, the usual approach is to apply an iterative solver. In such a solver it is critical
to compute the system operator $F(\mathcal{c})$ as efficiently as possible. At first glance, the tensors  $\mathcal{w}^{\mathbf{klm}},~\mathcal{f}^{\mathbf{klm}}$ 
could be of a large size and their direct computation appears to be intractable due to huge memory requirements. However, due to the finite support of B-spline functions these tensors have a large number of zeros, i.e.~they are sparse. Moreover, the non-zero values are localized around the grid nodes and are translation invariant. Therefore, they have a kernel-like structure, where  the width of these kernels depends on the B-spline degree.

Inside the domain, the kernels are translation invariant, and due to the B-spline separability property, these kernels are also separable:
\begin{eqnarray}
\hat{\mathcal{w}}_{\mathbf{klm}} = (1/h^2) \hat{\mathcal{w}}_{k_1 l_1 m_1} \cdots \hat{\mathcal{f}}_{k_d l_d m_d} + ... + 
(1/h^2) \hat{\mathcal{f}}_{k_1 l_1 m_1} \cdots \hat{\mathcal{w}}_{k_d l_d m_d} ; ~~~
\hat{\mathcal{f}}_{\mathbf{klm}} = \hat{\mathcal{f}}_{k_1 l_1 m_1} \cdots \hat{\mathcal{f}}_{k_d l_d m_d}.
\end{eqnarray}

\begin{figure}[H]
\scriptsize
\center
\def\svgscale{0.5}
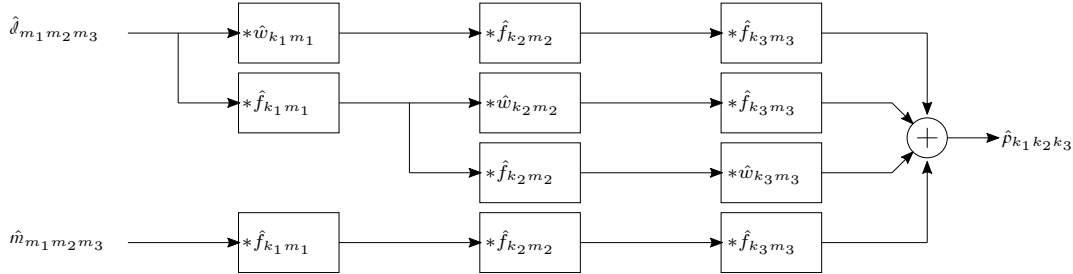
\caption{Domain computations are implemented as filtering algorithm} \label{figFiltAlg}
\normalsize
\end{figure}

Such a kernel-based decomposition structure of the system operator induces filtering-like algorithms which makes it possible to define algorithms from a signal processing viewpoint. The translation invariant kernels are the key component in these efficient filtering-like algorithms. Inside the domain, the operations are presented by multilinear convolution, with shift-invariant separable kernels, as depicted in Fig.~\ref{figFiltAlg}.

When the domain boundary is of arbitrary shape, the rectangular grid of the B-spline basis does not conform to the boundary geometry. Therefore, B-splines will be truncated on the boundary. These truncated B-splines result in non-separable kernels. We refer the reader to  \cite{Shulga2017, Shulga2018} where an efficient method for the integration of such non-separable kernels was proposed. Usually, the number of such kernels is small in comparison to the number of domain kernels. 

\subsection{Method Overview and High Performance Computational Strategies}

\begin{figure}[H]
\scriptsize
\def\svgscale{0.5}
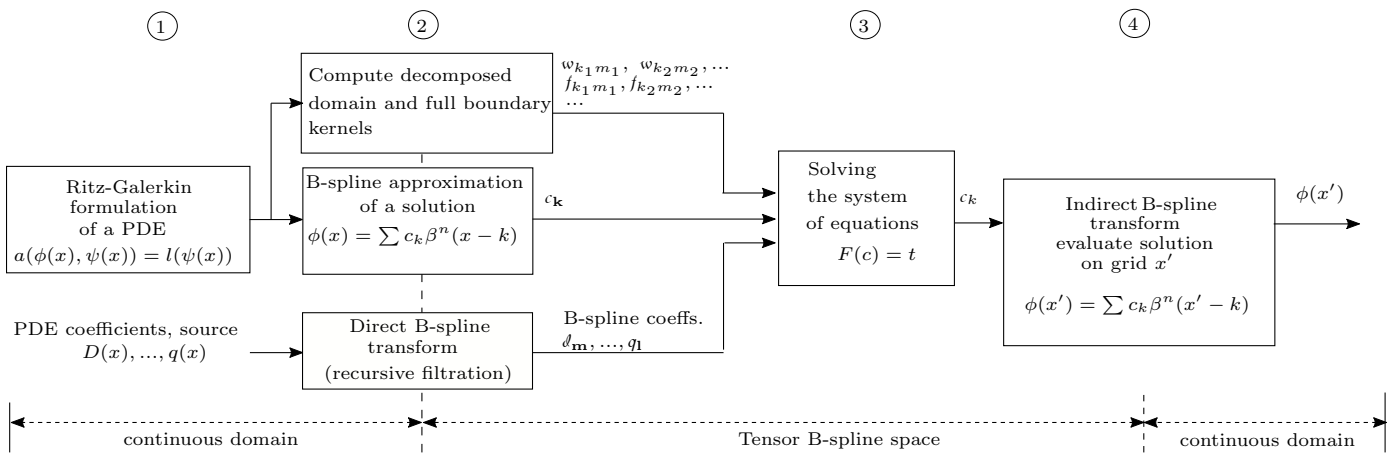
\caption{An overview of the Tensor B-spline numerical method for PDEs } \label{figGenTenSolv}
\normalsize
\end{figure}

An overview of the generic Tensor B-spline numerical method for  PDEs is given in Fig.~\ref{figGenTenSolv}. The method starts in the continuous domain with the Ritz-Galerkin formulation of a  PDE, given source and coefficients \raisebox{.5pt}{\textcircled{\raisebox{-.9pt} {1}}}. The continuous formulation is discretized with Tensor B-splines \raisebox{.5pt}{\textcircled{\raisebox{-.9pt} {2}}}. At this stage the domain and boundary kernels are computed and both the source and the coefficients are transformed into B-spline space via the direct B-spline transform. The obtained system of equations is solved using some suitable method. In this paper  we consider the conjugate gradient method targeted for large-scale problems \raisebox{.5pt}{\textcircled{\raisebox{-.9pt} {3}}}. The obtained coefficients of the solution are transformed back into the continuous domain by way of the indirect B-spline transform \raisebox{.5pt}{\textcircled{\raisebox{-.9pt} {4}}}.

\begin{figure}[H]
\centering     
\includegraphics[scale=0.25]{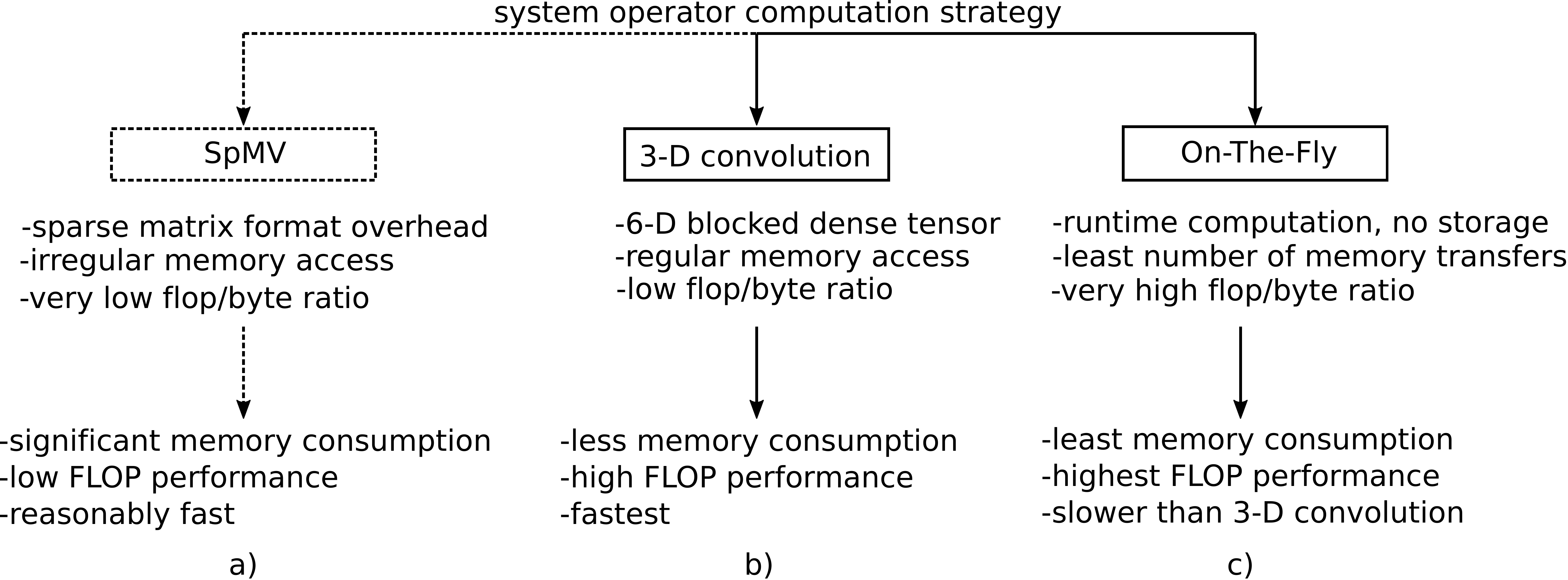}
\caption{Computation strategies of the Tensor B-spline method} \label{figStrategy}
\end{figure}

The system operator \eqref{tensForm} is the fundamental building block of iterative methods, and it defines  
the most time-consuming stage when solving large sparse linear systems. Fig.~\ref{figStrategy} depicts possible strategies provided by the Tensor B-spline method for this operation.

A standard sparse matrix-vector multiplication (SpMV) procedure (Fig.~\ref{figStrategy}, a) requires  assembling a sparse matrix from B-spline kernels. It was shown that this approach appears to be the least efficient \cite{Shulga2017, Shulga2018}. This follows from the overhead due to the sparse matrix format, from non-regular memory access, from a very low flop-to-byte ratio \cite{saule2013xeonphi, elafrou2017},  and from problems concerning load imbalance \cite{hou2017}. Since SpMV is a memory-bound procedure, performance optimizations do not overcome the issue of considerable memory consumption.
      
The use of tensor structures permits the implementation of more efficient computing algorithms. The first one uses a natural 6-D block tensor that reduces memory consumption and provides regular memory access. The 6-D block tensor $\mathcal{P}^{\mathbf{kl}} = 
\mathcal{d}_{\mathbf{m}} \mathcal{w}^{\mathbf{klm}} + 
\mathcal{m}_{\mathbf{m}} \mathcal{f}^{\mathbf{klm}} + 
(1/2) \mathcal{h}^{\mathbf{kl}} $ is assembled before the iterative solving. At each iteration it is multiplied with~$\mathcal{c}_{\mathbf{k}}$. This algorithm is the fastest \cite{Shulga2017,Shulga2018} but still is not feasible for large systems (Fig.~\ref{figStrategy}, b). The second, on-the-fly algorithm  computes the $ \mathcal{c}_{\mathbf{k}} \mathcal{d}_{\mathbf{m}} \mathcal{w}^{\mathbf{klm}} + 
\mathcal{c}_{\mathbf{k}} \mathcal{m}_{\mathbf{m}} \mathcal{f}^{\mathbf{klm}} + 
(1/2) \mathcal{c}_{\mathbf{k}} \mathcal{h}^{\mathbf{kl}} $ in the runtime and results in a significant reduction in  memory usage (Fig.~\ref{figStrategy}, c). It has a very high flop-to-byte ratio and gains from high floating point performance of CPUs and GPUs, as we are showing below.

\section{Method Performance} \label{secMethodPerf}

\begin{figure}[H]
\centering     
\includegraphics[scale=0.43]{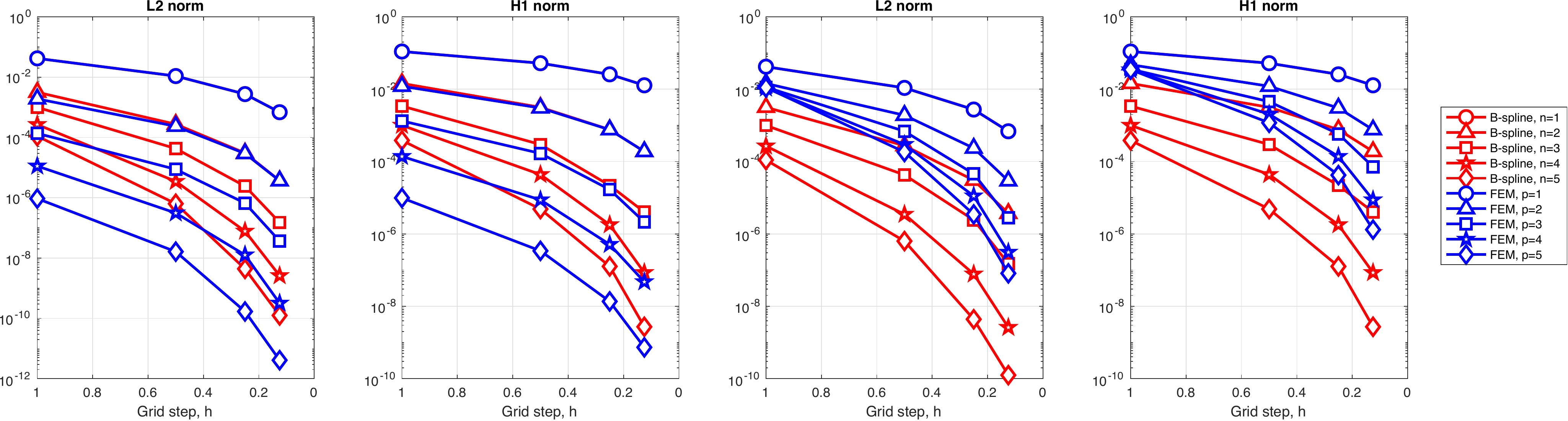}
\begin{picture}(400,10)
\put(-0,0){a)}
\put(115,0){b)}
\put(240,0){c)}
\put(360,0){d)}
\end{picture}
\caption{$L^{2}$ and $H^{1}$ norms of errors between the reference solution and B-spline and the reference and FEM solution. FEM uses additional nodes (a), (b). FEM uses the same number of nodes as B-splines use (c), (d).} \label{L2H1_1d}
\end{figure}
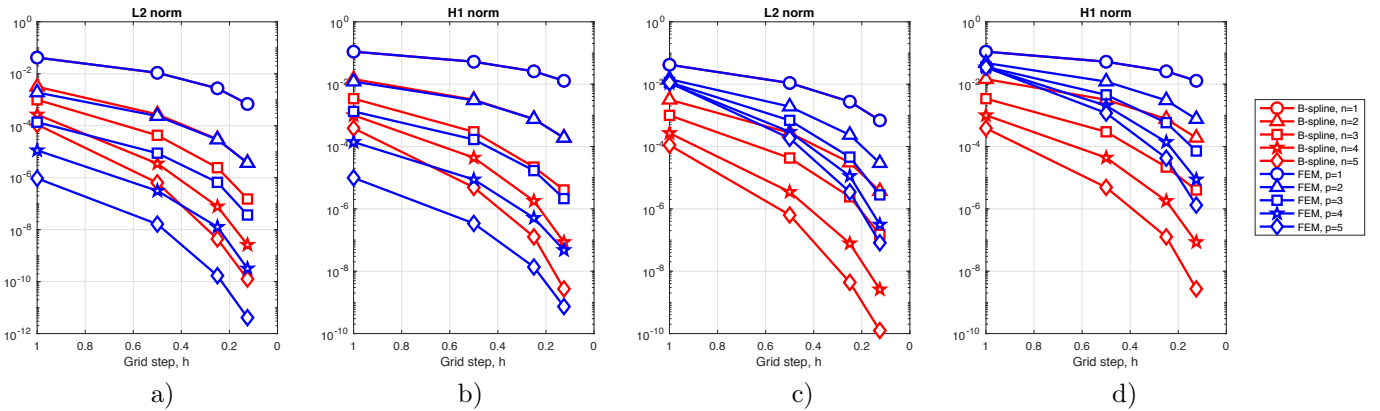

We begin with a didactic one-dimensional example that compares the accuracies of a Tensor B-spline solver and an FEM solver. For a comprehensive comparison with FEM in 2-D and 3-D, we refer the reader to \cite{Shulga2017,Shulga2018}. The problem is defined by a Diffusion PDE with Robin BC:
\begin{eqnarray}
- \nabla \cdot (D(\mathbf{x})\nabla \varphi(\mathbf{x})) + \mu_{a}(\mathbf{x}) \varphi(\mathbf{x}) &=& q(\mathbf{x}), ~\mathbf{x} \in ~ \Omega  \\
2D(\mathbf{x}) (\nabla \varphi(\mathbf{x}) \cdot \mathbf{n}) + \varphi(\mathbf{x}) &=& 0, ~\mathbf{x} \in \partial \Omega.
\end{eqnarray}
The domain limits are $[-25,25]$, the source is $q(x) = \exp(-(\frac{x}{2})^2)$. The grid step $h$ was decreased using the expression $h(\mu)=2^{-\mu},~\mu=\{0, 1, 2, 3\}$. Note that both method's bases correspond to the ones shown in Fig.\  \ref{figFEMSplineBasis}. The $L^2$ and $H^2$ errors (defined in Appendix I) are computed between the numerical solutions and the analytic reference solution. Two situations were studied, in which: 1) additional nodes where introduced in FEM, and 2) FEM was forced to have the same number of nodes as the Tensor B-spline method. The plots of errors are shown in Fig.~\ref{L2H1_1d}.

\scriptsize
\begin{table}[H]
\centering
\small
\caption{Number of operations for  Fig.~\ref{L2H1_1d} }  \label{tblFlop1d}
\begin{tabular}{ r    c c c c c | c c c c c | c c c c c }

& \multicolumn{5}{ c|}{B-spline} & \multicolumn{5}{c|}{ FEM (with additional nodes)}  & \multicolumn{5}{c}{ FEM (same number of nodes)}  \\
\hline
h        & n=1 & n=2 & n=3 & n=4 & n=5              & p=1 & p=2 & p=3 & p=4 & p=5                  & p=1 & p=2 & p=3 & p=4 & p=5  \\
\hline
1        & 251   & 465  & 665   & 895 &  1095      & 251  & 701   &  1351  & 2201  &3251            &251 &351  & 433  & 529  & 651\\
0.5    & 501   & 915   & 1315   & 1745 &  2145   & 501   & 1401 & 2701 & 4401   & 6501          &501  & 701 & 892  & 1101 & 1301\\ 
0.25  & 1001  & 1815 & 2615  & 3445 &  4245  & 1001 & 2801 & 5401 & 8801 & 13001          &1001 & 1401 & 1783 & 2201 & 2601\\
0.125 & 2001 & 3615 & 5215 & 6845 &  8445  & 2001 & 5601 & 10801 & 17601 & 26001       &2001 & 2801 & 3592 & 4401& 5201\\
\hline
\end{tabular}
\end{table}
\normalsize

From Fig.~\ref{L2H1_1d} one can observe that a high-order FEM with additional nodes is more accurate than a high-order Tensor B-spline, but at the same time, the number of operations for  a high-order FEM  grows dramatically (Table~\ref{tblFlop1d}). With the same number of nodes, a high-order Tensor B-spline is more accurate while requiring only slightly more operations.

\begin{figure}[H]
\centering     
\includegraphics[scale=0.5]{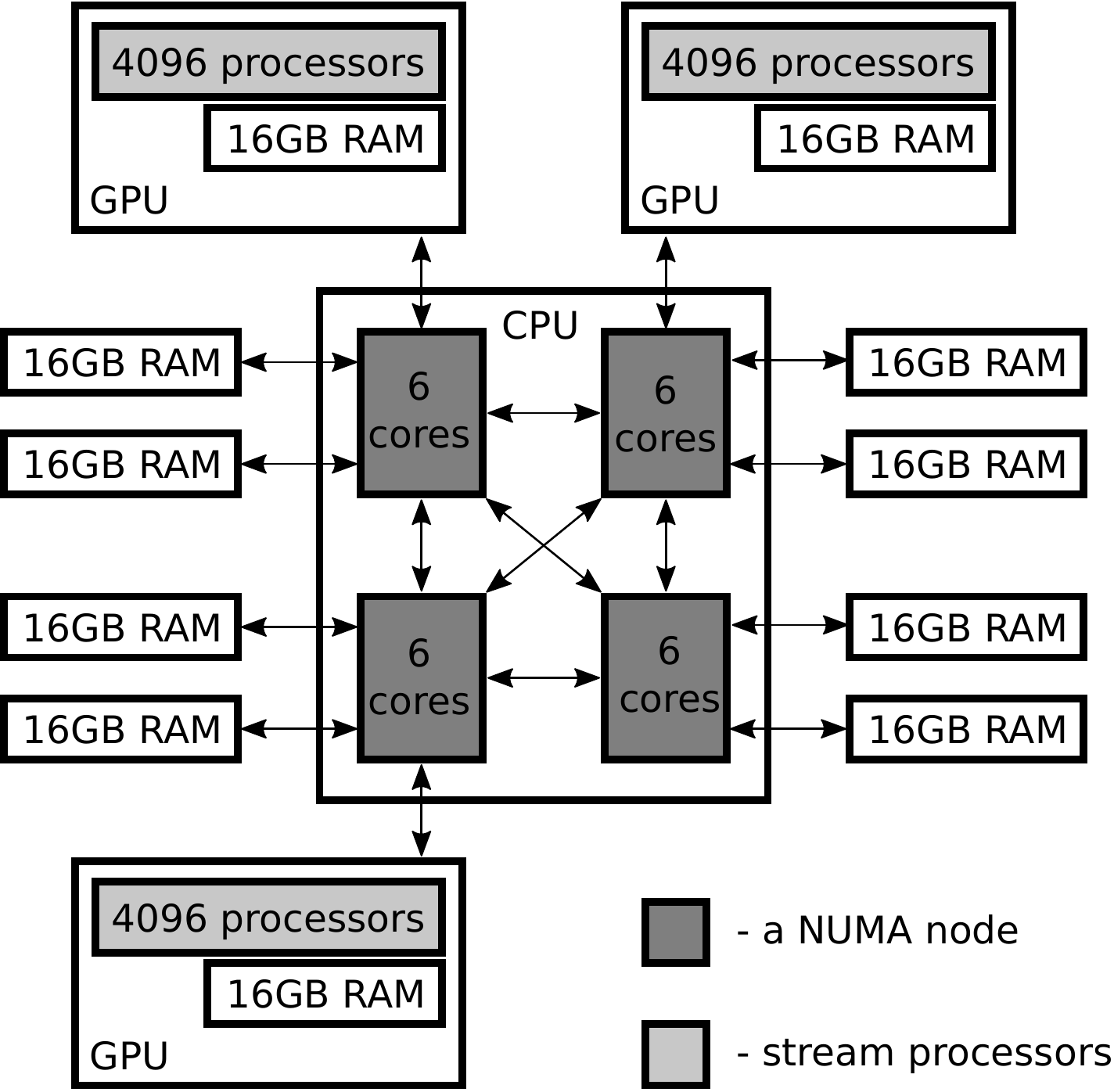}
\caption{The target high-performance workstation consisting of the AMD EPYC 7401P 2.0 GHz 24 cores CPU, 128GB RAM, and three AMD Radeon Vega Frontier Edition 16 GB RAM GPUs.} \label{figEpyc}
\end{figure}

The next benchmark considers a three-dimensional case. The performance of the system operator, stated to be the most critical part of the Tensor B-spline solver, was estimated on a heterogeneous workstation. A simplified diagram of the architecture of the workstation is shown in Fig.~\ref{figEpyc}. The workstation's CPU AMD EPYC 7401P has four non-uniform memory access (NUMA) nodes (shown in dark gray). Each node is connected to its own random-access memory (RAM) domain and to the other nodes. There are 128 GB ($4 \times 2 \times 16 $ GB) of RAM in total. Each node contains a multi-core processor with floating-point SIMD units and several levels of cache memory. Three nodes are connected to GPUs AMD Radeon Vega Frontier Edition. Each GPU has 4096 stream processors (shown in light gray) and 16 GB of RAM.

\begin{figure}[H]
\centering     
\includegraphics[scale=0.4]{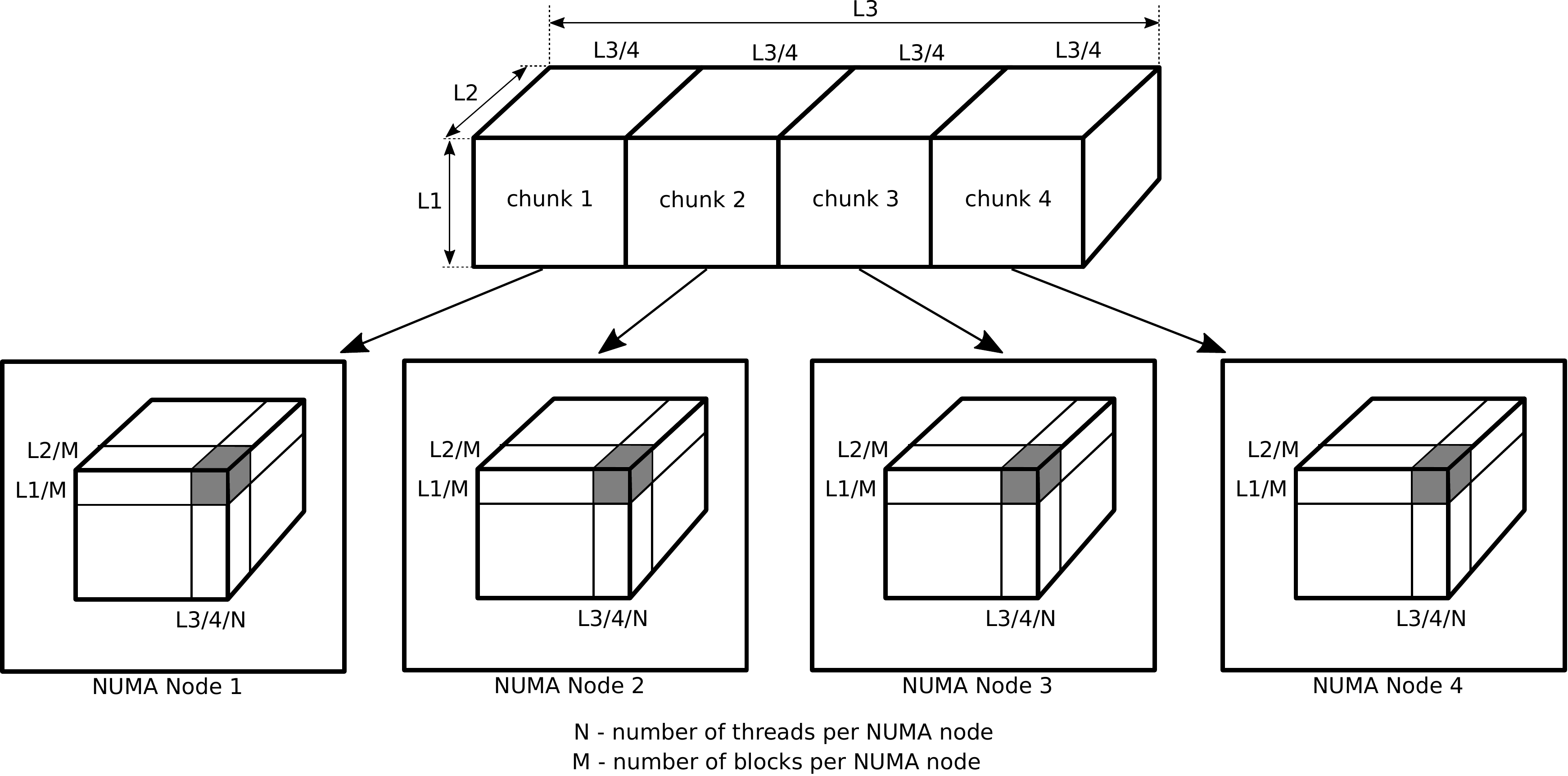}
\caption{The pattern of data distribution for CPU-based computations} \label{figDataDist}
\end{figure}

The platform's massive parallelism and non-uniform memory access could challenge an efficient implementation of a numerical algorithm. However, the Tensor B-spline method allows us to use data parallelism in a very intuitive way.  Fig.~\ref{figDataDist} depicts a possible pattern of data distribution that was applied to the  CPU-based computations in this paper. The rectangular domain of size $L1 \times L2 \times L3$ is divided into rectangular sub-domains (chunks) of size  $L1 \times L2 \times L3/4$. The chunks are bound to NUMA nodes and further divided into blocks of size $L1/M \times L2/M \times L3/4/N$ in order to be processed independently and simultaneously by threads involving SIMD instructions. One can adjust the number of threads  per node ($N$) as well as the number of blocks per node ($M$) to tune the performance and CPU load. The  majority of computations is represented by filtration of input data (see Fig.~\ref{figFiltAlg}) performed with the use of fused multiply-add (FMA) SIMD instructions.

The workstation runs 64-bit Linux Ubuntu 16.04.5 LTS. We used GCC 8.1.0 compiler, POSIX threads for parallelization, and libnuma library for memory management and threads binding to NUMA nodes, AVX instructions for vectorized floating point computations. The GPU-based computations were performed using an Active Oberon tensor runtime with OpenCL support. 

The important property of any numerical algorithm is its scalability, since it allows us to understand the potential performance that can be achieved on different platforms and even supercomputers. The scalability of the method is presented in Table~\ref{tblScalability}, showing the  performance of the on-the-fly algorithm for different numbers of  exploited NUMA nodes and threads. 

\begin{table}[H]
\centering
   \caption{On-the-fly algorithm scalability performance  on the AMD AMD EPYC 7401P CPU.  Grid of size 240x240x960, cubic splines.} \label{tblScalability}
\def\arraystretch{1}%
\begin{tabular}{c c c c }

\hline
NUMA nodes            & threads				& time, s     & GFLOPS      \\
\hline						 
{\multirow{3}{*}{1} } & 4                         & 84.93    & 36                \\
                                  & 8                         & 72.43    & 42                \\
                                  & 12                       & 50.90    & 60                \\
\hline 
{\multirow{3}{*}{2} }  & 8                         & 48.10    & 64               \\
                                   & 16                       & 31.67    & 96                \\
                                   & 24                       & 26.00    & 117              \\
\hline 
{\multirow{3}{*}{4} }   & 12                        & 27.78    & 110             \\
                                    & 24                       & 17.55    & 174              \\
                                     & 48                       & 14.69    & 208             \\
\hline 
\end{tabular}
\end{table}

The Tensor B-spline method shows a significant increase in  performance when more and more NUMA nodes and threads are used. Table~\ref{tblPerformance} shows the performance of the on-the-fly algorithm for different B-spline degrees, grids and precisions. 

\begin{table}[H]
\centering
   \caption{On-the-fly algorithm performance on the AMD EPYC CPU only.} \label{tblPerformance}
\def\arraystretch{1}%
\begin{tabular}{ c|c|r| c r r| c r r }

\multicolumn{3}{ c|}{} & \multicolumn{3}{|c|}{ double precision}   & \multicolumn{3}{|c }{ single precision}  \\

\hline
n    & Grid(L1xL2xL3)						 & GFLOP     & MEM R/W, GB     &  time, s    &  GFLOPS    & MEM R/W, GB     &  time, s    &  GFLOPS   \\
\hline						 
1   &   240x240x960   & 97            &   2.11/0.42          &   1.07        &       91         &   1.06/0.21           & 0.87        &  111                  \\
2   &  240x240x960 	 & 781           & 2.16/0.44           &   6.93       &      113        &    1.08/0.22          & 3.65        &   213            \\
3   &  240x240x960   & 3054        & 2.22/0.44          &  15.11       &      202        &    1.11/0.22           &  9.87        &   309            \\ 

\hline 
1   &   480x480x960   & 388          &   8.38/1.68        &   3.81       &     102         &    4.19/0.84        &   3.72          &  104             \\
2   &  480x480x960 	 &  3124        & 8.52/1.70           & 21.75       &     144        &     4.26/0.85        & 14.14         &   221                 \\
3   &  480x480x960   &  12214      & 8.66/1.73          &  59.81      &      204       &     4.33/0.87         & 40.10       &    305                      \\ 

\hline 
1   &   960x960x960   &  1552       & 33.37/6.67       & 14.34        &     108        &     16.69/3.34       & 13.36       & 116                         \\
2   &  960x960x960 	 &   12496    & 33.79/6.76       & 83.85        &    149         &     16.89/3.38      &  57.05      & 219                         \\
3   &  960x960x960   &   48857   & 34.21/6.84        &237.90       &    205        &      17.10/3.42      &  163.76     &   298                       \\ 

\hline 
1   &   1200x1200x1200   &  3031       & 65.02/13.00   &  28.70    & 106            &     32.51/6.50       & 26.62       & 114                         \\
2   &  1200x1200x1200   &   24406    & 65.67/13.13    &  172.94     & 141             &     32.83/6.57      &  111.65      & 219                         \\
3   &  1200x1200x1200   &   95424    & 66.32/13.26   &   474.72    &  201          &      33.16/6.63     & 324.12      & 294                         \\ 

\hline 
\end{tabular}
\end{table}

Table~\ref{tblPerformance} shows that Tensor B-spline method achieved  200 GFLOPS in double precision and 300 GFLOPS in single precision when cubic basis functions are used. Despite of the very large problem sizes, the on-the-fly strategy provides conservative memory usage, thus the problem of 1.7 billion nodes requires only  66.3 GB to read and 13.26 GB to write.  The amount of used memory almost does not depend on the B-spline degree.

Table~\ref{tblGpu} depicts GPU performance of the on-the-fly algorithm for single precision.

\begin{table}[H]
\centering
   \caption{AMD RADEON Vega GPU  performance for the on-the-fly algorithm} \label{tblGpu}
\def\arraystretch{1}%
\begin{tabular}{ c| c | c | cc | cc | cc  }

\multicolumn{3}{ c|}{}       & \multicolumn{2}{c|}{One GPU} & \multicolumn{2}{c|}{Two GPUs} & \multicolumn{2}{c}{Three GPUs}  \\
\hline
n & Grid						 & Memory,GB&  time               &      GFLOPS     &  time               &      GFLOPS         &  time               &      GFLOPS              \\
\hline						 
3 & 72x72x72		&     0.00695   &   68 ms           &       212       &   71 ms             &        204               &        -                  &      -                                  \\
3 &  144x144x144 	&     0.0556     &   510 ms            &      227      &   512 ms         &      226                 &        -                  &       -                               \\
3 &  258x258x258   &    0.3198     &    3.821 s         &      243      &    1.912 s         &      485                 &   1.333 s             &       711                       \\ 
3 & 522x522x522   &    2.649             &   21.079 s        &      352     &      10.54 s       &      703                  &   7.4486 s          &       1055                    \\
\hline
\end{tabular}

\end{table}

In the last example, we present the solution of the heat equation in steady state with Mixed BC (Dirichlet and Neumann)
\begin{eqnarray}
- \nabla \cdot (k(\mathbf{x})  \nabla u(\textbf{x})) &=& g(\mathbf{x}) , \quad \in \Omega  \\
u(\textbf{x}) &=& u_0, \quad \in \partial \Omega_1 \\
\nabla \varphi(\mathbf{x}) \cdot \mathbf{n} &=& 0, \quad \in \partial \Omega_2,
\end{eqnarray}
where $u(\textbf{x})$ is the temperature, $k(\mathbf{x}) $ is the thermal conductivity, and $g(\mathbf{x}) $ is the rate of heat generation. 

\begin{figure}[H]
\centering     
\includegraphics[scale=1.5,angle=-90,trim={6cm 0cm 6cm 0cm},clip]{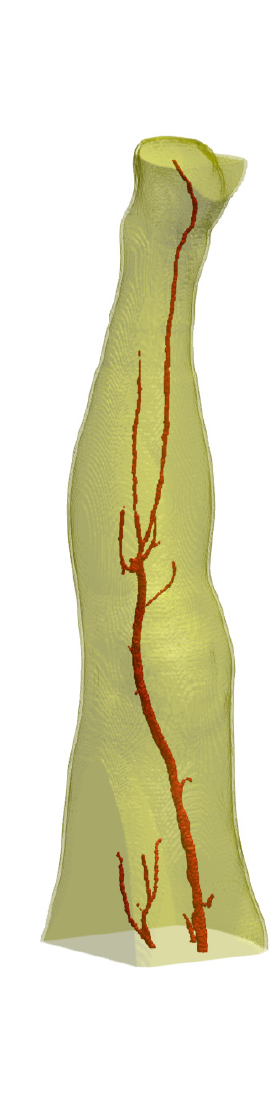}
\caption{The volumetric visualization of the leg's skin and arteria obtained from a CT dataset} \label{figLegDomain}
\end{figure}

The domain of numerical computations is obtained by segmentation of a computed tomography (CT) scan of a human leg. Skin, blood vessels, surrounding tissue, and the internal volume of the leg were separated based on gray-scale values, holes were closed and air is added around the skin. The resulting domain size is $600 \times 600 \times 2400$ ($\approx 0.8$ billions of unknowns). The grid step (resolution) is $h=0.3$ mm. The visualization of the leg's skin and arteria is shown in Fig.~\ref{figLegDomain}.

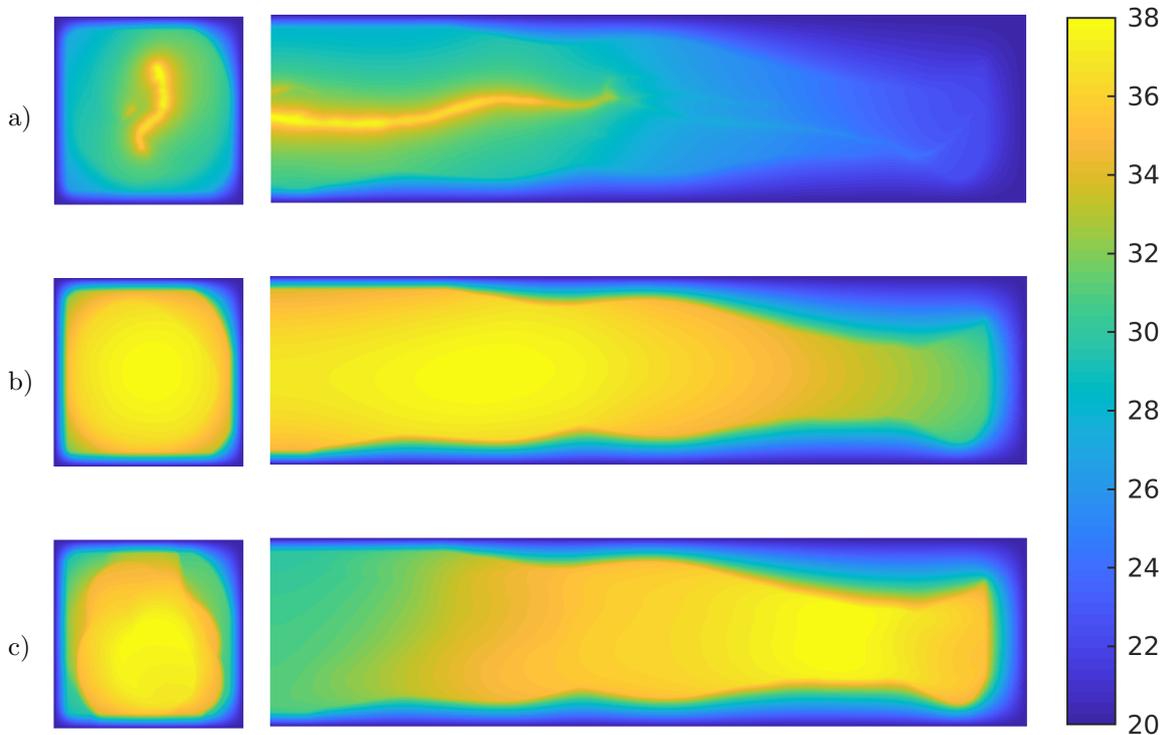
\begin{figure}[H]
\centering     
\includegraphics[scale=1.2]{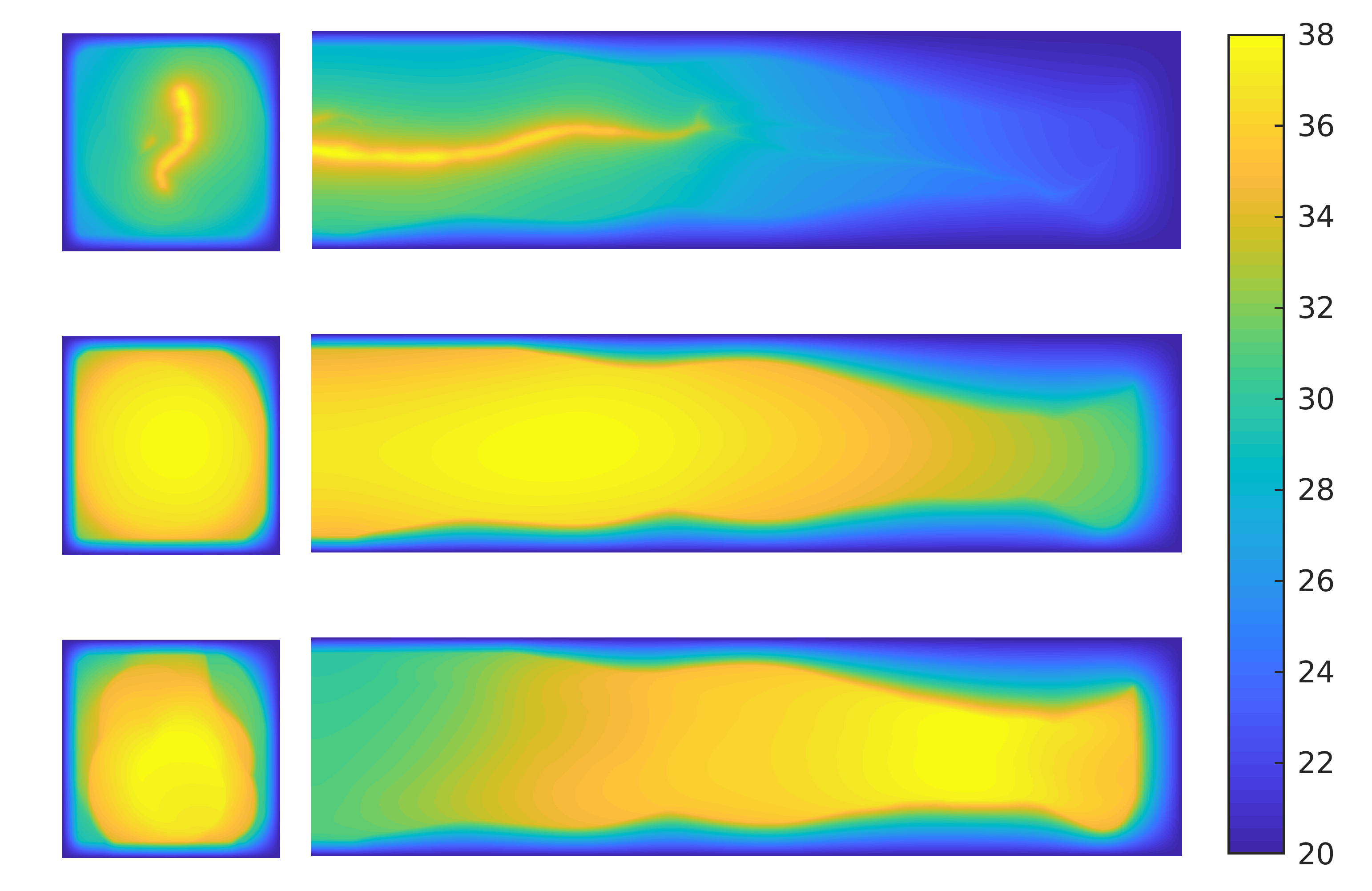}
\begin{picture}(400,10)
\put(-20,50){c)}
\put(-20,150){b)}
\put(-20,250){a)}
\end{picture}
\caption{The solution of the heat equation on leg domain. The color bar shows the temperature values in $\degree C$. The source of heat is a) arteria, b) tissue, c) skin.} \label{figLegSol}
\end{figure}

For the numerical computations, a data distribution pattern was used as depicted in Fig.~\ref{figDataDist}. There are  $98.5~\%$ of domain kernels and $1.5~\%$ of boundary kernels. A linear B-spline basis was used. The solution was computed in double precision.  The Dirichlet BC was approximated with Robin BC using the boundary penalty method \cite{barrett1986}. The Dirichlet BC sets the temperature on the boundary to $20~\degree C$ .  The advantage of this method is that it does not require the basis functions to vanish on the boundary. The solver uses the parallel implementation of the on-the-fly algorithm of the system operator, because  neither sparse matrices nor 6-D tensors would fit into the 128 GB of RAM available. A parallel version of the conjugate gradient algorithm with Jacobi preconditioner was used. In order to speed-up the convergence on such a large grid, the solver was initialized with a solution obtained on a reduced grid of size $120 \times 120 \times 480$. The solution is presented in Fig.~\ref{figLegSol} as maximum-projections images.  Each solution was computed in 4 h 4 min for 1000  conjugate gradient iterations.

\section{Discussion} \label{secDiscussion}

The numerical solution of PDEs is a computationally intensive task in many scientific and engineering applications. A general requirement for numerical PDE solvers is the accurate representation of the underlying continuous mathematical model, and from a practical point of view, it is absolutely crucial to utilize the computational resources available as efficiently as possible. Both aspects can only be addressed successfully, if the solution strategy used is flexible enough to adapt to the intrinsic problem structure, and if the algorithms used can adapted easily to the architecture of the hardware. In this paper we studied the Tensor B-spline method which for many practical problems indeed meets these requirements, thereby providing a general and flexible numerical PDE solver. Its key properties are (i) an accurate and flexible discretization on coefficients, sources and the solution, (ii) a highly efficient integration strategy, (iii) the availability of efficient algorithms. 

\begin{enumerate}

\item \textbf{Domain Discretization}

Any element-based numerical method relies on a specific domain discretization. This discretization is represented either by elements (in case of FEM) or by grid cells (for Tensor B-splines). A key aspect of most FEM-based approaches is the use of a domain-conforming mesh (which is often unstructured). In FEM, meshing automation is a major field of study, and finding suitable automation strategies is highly challenging in practice \cite{bathe2019}. Meshless methods, on the other hand, require much less effort when it comes to discretizing the domain, but the process of numerical integration  typically involves very high computational costs \cite{De2000}. Contrary to these approaches, tensor B-splines on regular grids have the advantage of regularity and shift-invariance in the domain. For these methods,  however, boundary integration is a challenging problem, as well as the use of Dirichlet BC on non-rectangular domains. In  \cite{Shulga2017,Shulga2018} an efficient integration method was proposed, based on the Divergence Theorem. Interestingly, recent advantages in FEM also consider similar approaches, where an easy-to-build regular mesh is used inside the domain, and the conforming mesh is adjusted to the domain boundaries \cite{bathe2019}. Due to intrinsic  B-spline multi-resolution properties, however, local grid refinements  -- similar to such mesh adaptations in FEM -- are also possible in Tensor B-spline approaches. 

\item \textbf{Data representation: Tensors and Matrices}

Tensors represent multidimensional data in a natural way in comparison to matrices, where the initial spatial structure of the data is folded  \cite{Cichocki2015, Sidiropoulos2017, Morozov2010}. This important aspect of tensors allows one to identify repeated patterns, local kernels  and convolution procedures in order to derive efficient data processing algorithms. The index-based notation for tensors can suggest an optimal order of tensor operations in a natural way. The kernel-based tensor array format avoids the overhead of the sparse matrix format, where indices of the non-zero elements have to be stored and accessed. Moreover,  SpMV operations are known for their low computational efficiency \cite{saule2013xeonphi, elafrou2017}. Tensor B-spline methods allow us to use tensor-based computational kernels (3-D convolutions and on-the-fly) that are more efficient than SpMV \cite{Shulga2017, Shulga2018}. 

\item \textbf{High-Order Basis}

The use of high-order approximations is beneficial, because faster convergence rates can be achieved \cite{Guo1986,Mitchell2014,Mitchell2015}.  For a specified level of accuracy, high-order methods are typically more efficient than low-order methods. Even more important, high-order Tensor B-splines are usually more efficient than high-order polynomial FEM basis, see \cite{Shulga2017,Shulga2018}. Further, such high-order B-splines require less memory and operations per iteration, and usually the convergence rate is faster, too.

\item \textbf{High-Performance Computing and Applications}

Some benefits of the use of tensor B-splines for high performance computing have already been shown in \cite{Shulga2017, Shulga2018}. In this work, we present a parallel high-performance implementation for multi-core CPU with shared memory and GPUs, and we show numerical results for a large-scale problem consisting of 0.8 billion nodes. We demonstrated that tensor B-spline methods are indeed capable of solving such large-scale problems within reasonable time limits, whereas standard FEM methods run into severe memory problems. Recent advances in element-based methods for large-scale weather forecasting problems \cite{marras2016} suggest a huge potential application field for such tensor B-spline methods.  

\item \textbf{Dirichlet BC}

The application of the Dirichlet BC is challenging in methods with basis functions that do not conform to the domain boundary. As was proposed in \cite{Hollig2003}, one can use a weighted basis that vanishes on the boundary. In this work we showed that tensor B-spline solvers can successfully make use of the boundary penalty method \cite{barrett1986}. In this method the Dirichlet BC is approximated with Robin BC weighted by some penalty factor.

\item \textbf{Domain-Specific Hardware}

Despite ubiquitous utilization of general-purpose CPUs and GPUs, recent studies suggest domain-specific hardware as a future solution for efficient computations \cite{Jouppi2017}. An example of this kind is Google's Tensor Processing Unit containing 65,536 8-bit MAC matrix multiply units with a peak throughput of 92 TeraOps/second. The specialized computing units of a such hardware rely on specific properties of the data structure and data flow, as well as on specific data locality patterns \cite{Qadeer2013}.  One interesting aspect of tensor B-spline methods is that they can be naturally deployed on such domain-specific hardware, a combination which -- at least in our view -- has a huge potential in future high-performance computing problems. On the other hand, tensor B-spline methods allow us to solve large-scale problems solution using limited compute resources. This might also suggests the use of these methods in embedded systems with their typical restrictions in terms of on-chip memory and flops-per-watt ratio.

\end{enumerate}

\section{Conclusion} \label{secConclusion}

In this paper the tensor B-spline method was described, with specific emphasis on its mathematical foundations and implications for developing high-performance algorithms. There is a clear trend that over the years both cutting-edge as well as routine PDE simulations require more and more compute resources. Classical FEM approaches to such simulations, however, reside in a relatively inflexible framework of mesh-based polynomial bases and sparse matrix computations. While in recent years, many research efforts addressed the problem of mesh optimization and SpMV computations, we presented a conceptually different approach: the combination of tensor algebra and B-spline bases allows for a mesh-free, high-order, parallel, matrix-free numerical method for solving PDEs. Compared to FEM methods, tensor B-splines share several advantages, both on the conceptual side and on the application side: the mathematical elegance of computational tensor algebra leads to flexible and transparent models, that respect the intrinsic multi-dimensional structure of the problem in a very natural way. Further, Tensor B-splines offer accurate numerical solutions at significantly reduced computational costs and memory requirements. This latter advantage makes these methods particularly interesting for resource-constraint applications in, for instance, embedded systems. 

\section{Acknowledgment}

The authors thank Prof. Christoph J. Zech from University Hospital of Basel for providing CT data sets used for numerical simulations presented in this paper.

\section{Appendix I}

The $L^2-$norm for a function $f(x) \in  \Omega$ is defined as 
\begin{equation}
||f(x)||_{\mathbb{L}^{2}(\Omega)} = \sqrt{\int_{\Omega} |f(x)|^2 dx}
\end{equation}

The $H^{1}-$norm (energy norm) for a function $f(x) \in  \Omega$ is defined as
\begin{equation}
||f(x)||_{\mathbb{H}^{1}(\Omega)} = \sqrt{ \int_{\Omega} |f(x)|^2 dx + \int_{\Omega} | \nabla f(x) |^2 dx}
\end{equation}

\end{document}

%% file: Domain.pdf_tex
\begingroup%
  \makeatletter%
  \providecommand\color[2][]{%
    \errmessage{(Inkscape) Color is used for the text in Inkscape, but the package 'color.sty' is not loaded}%
    \renewcommand\color[2][]{}%
  }%
  \providecommand\transparent[1]{%
    \errmessage{(Inkscape) Transparency is used (non-zero) for the text in Inkscape, but the package 'transparent.sty' is not loaded}%
    \renewcommand\transparent[1]{}%
  }%
  \providecommand\rotatebox[2]{#2}%
  \ifx\svgwidth\undefined%
    \setlength{\unitlength}{196.25986495bp}%
    \ifx\svgscale\undefined%
      \relax%
    \else%
      \setlength{\unitlength}{\unitlength * \real{\svgscale}}%
    \fi%
  \else%
    \setlength{\unitlength}{\svgwidth}%
  \fi%
  \global\let\svgwidth\undefined%
  \global\let\svgscale\undefined%
  \makeatother%
  \begin{picture}(1,0.8185207)%
    \put(0,0){\includegraphics[width=\unitlength,page=1]{Domain.pdf}}%
    \put(0.28864343,0.76045236){\color[rgb]{0,0,0}\makebox(0,0)[lb]{\smash{$\partial \Omega$}}}%
    \put(0.43877246,0.30733611){\color[rgb]{0,0,0}\makebox(0,0)[lb]{\smash{$\Omega$}}}%
    \put(0,0){\includegraphics[width=\unitlength,page=2]{Domain.pdf}}%
    \put(0.01841142,0.05621142){\color[rgb]{0,0,0}\makebox(0,0)[lb]{\smash{$\mathbf{n}$}}}%
  \end{picture}%
\endgroup%

%% file: SplineProc.pdf_tex
\begingroup%
  \makeatletter%
  \providecommand\color[2][]{%
    \errmessage{(Inkscape) Color is used for the text in Inkscape, but the package 'color.sty' is not loaded}%
    \renewcommand\color[2][]{}%
  }%
  \providecommand\transparent[1]{%
    \errmessage{(Inkscape) Transparency is used (non-zero) for the text in Inkscape, but the package 'transparent.sty' is not loaded}%
    \renewcommand\transparent[1]{}%
  }%
  \providecommand\rotatebox[2]{#2}%
  \ifx\svgwidth\undefined%
    \setlength{\unitlength}{842.51151973bp}%
    \ifx\svgscale\undefined%
      \relax%
    \else%
      \setlength{\unitlength}{\unitlength * \real{\svgscale}}%
    \fi%
  \else%
    \setlength{\unitlength}{\svgwidth}%
  \fi%
  \global\let\svgwidth\undefined%
  \global\let\svgscale\undefined%
  \makeatother%
  \begin{picture}(1,0.14718125)%
    \put(0,0){\includegraphics[width=\unitlength,page=1]{SplineProc.pdf}}%
    \put(0.01112605,0.11440636){\color[rgb]{0,0,0}\makebox(0,0)[lb]{\smash{Input signal $s(x)$}}}%
    \put(0,0){\includegraphics[width=\unitlength,page=2]{SplineProc.pdf}}%
    \put(0.22350328,0.12938172){\color[rgb]{0,0,0}\makebox(0,0)[lb]{\smash{Direct B-spline}}}%
    \put(0.22554537,0.10964153){\color[rgb]{0,0,0}\makebox(0,0)[lb]{\smash{transform}}}%
    \put(0.23099094,0.08853995){\color[rgb]{0,0,0}\makebox(0,0)[lb]{\smash{(recursive}}}%
    \put(0.23235233,0.0681191){\color[rgb]{0,0,0}\makebox(0,0)[lb]{\smash{filtration)}}}%
    \put(0,0){\includegraphics[width=\unitlength,page=3]{SplineProc.pdf}}%
    \put(0.39027384,0.12870102){\color[rgb]{0,0,0}\makebox(0,0)[lb]{\smash{Process signal }}}%
    \put(0.39027384,0.10828014){\color[rgb]{0,0,0}\makebox(0,0)[lb]{\smash{in space of B-splines}}}%
    \put(0.3916352,0.08241368){\color[rgb]{0,0,0}\makebox(0,0)[lb]{\smash{(manipulate with  }}}%
    \put(0.39231593,0.05858932){\color[rgb]{0,0,0}\makebox(0,0)[lb]{\smash{coefficients $c_{k}$ )  }}}%
    \put(0,0){\includegraphics[width=\unitlength,page=4]{SplineProc.pdf}}%
    \put(0.61518933,0.12938172){\color[rgb]{0,0,0}\makebox(0,0)[lb]{\smash{Indirect B-spline}}}%
    \put(0.61723142,0.10964153){\color[rgb]{0,0,0}\makebox(0,0)[lb]{\smash{transform}}}%
    \put(0.622677,0.08853995){\color[rgb]{0,0,0}\makebox(0,0)[lb]{\smash{(convolution with}}}%
    \put(0.62403839,0.0681191){\color[rgb]{0,0,0}\makebox(0,0)[lb]{\smash{sampled B-spline)}}}%
    \put(0,0){\includegraphics[width=\unitlength,page=5]{SplineProc.pdf}}%
    \put(0.82120398,0.11440638){\color[rgb]{0,0,0}\makebox(0,0)[lb]{\smash{Output signal $s'(x)$}}}%
    \put(0,0){\includegraphics[width=\unitlength,page=6]{SplineProc.pdf}}%
    \put(0.04243809,0.01036297){\color[rgb]{0,0,0}\makebox(0,0)[lb]{\smash{continuous domain}}}%
    \put(0.19642228,-1.4427789){\color[rgb]{0,0,0}\makebox(0,0)[lt]{\begin{minipage}{0.2355163\unitlength}\raggedright \end{minipage}}}%
    \put(0.23006746,-1.51006927){\color[rgb]{0,0,0}\makebox(0,0)[lt]{\begin{minipage}{0.20187112\unitlength}\raggedright \end{minipage}}}%
    \put(0.42892843,0.00959629){\color[rgb]{0,0,0}\makebox(0,0)[lb]{\smash{discrete domain}}}%
    \put(0,0){\includegraphics[width=\unitlength,page=7]{SplineProc.pdf}}%
    \put(0.76827601,0.01325092){\color[rgb]{0,0,0}\makebox(0,0)[lb]{\smash{continuous domain}}}%
    \put(0,0){\includegraphics[width=\unitlength,page=8]{SplineProc.pdf}}%
  \end{picture}%
\endgroup%

%% file: Tensor.pdf_tex
\begingroup%
  \makeatletter%
  \providecommand\color[2][]{%
    \errmessage{(Inkscape) Color is used for the text in Inkscape, but the package 'color.sty' is not loaded}%
    \renewcommand\color[2][]{}%
  }%
  \providecommand\transparent[1]{%
    \errmessage{(Inkscape) Transparency is used (non-zero) for the text in Inkscape, but the package 'transparent.sty' is not loaded}%
    \renewcommand\transparent[1]{}%
  }%
  \providecommand\rotatebox[2]{#2}%
  \ifx\svgwidth\undefined%
    \setlength{\unitlength}{443.24999884bp}%
    \ifx\svgscale\undefined%
      \relax%
    \else%
      \setlength{\unitlength}{\unitlength * \real{\svgscale}}%
    \fi%
  \else%
    \setlength{\unitlength}{\svgwidth}%
  \fi%
  \global\let\svgwidth\undefined%
  \global\let\svgscale\undefined%
  \makeatother%
  \begin{picture}(1,0.34852017)%
    \put(0,0){\includegraphics[width=\unitlength,page=1]{Tensor.pdf}}%
    \put(0.01490636,0.26364031){\color[rgb]{0,0,0}\makebox(0,0)[lb]{\smash{scalar}}}%
    \put(0,0){\includegraphics[width=\unitlength,page=2]{Tensor.pdf}}%
    \put(0.08344176,0.26314311){\color[rgb]{0,0,0}\makebox(0,0)[lb]{\smash{vector}}}%
    \put(0.18801483,0.26219244){\color[rgb]{0,0,0}\makebox(0,0)[lb]{\smash{matrix}}}%
    \put(0.33877644,0.10860686){\color[rgb]{0,0,0}\makebox(0,0)[lb]{\smash{3-way}}}%
    \put(0,0){\includegraphics[width=\unitlength,page=3]{Tensor.pdf}}%
    \put(0.59040231,0.11087405){\color[rgb]{0,0,0}\makebox(0,0)[lb]{\smash{4-way}}}%
    \put(0,0){\includegraphics[width=\unitlength,page=4]{Tensor.pdf}}%
    \put(0.01212004,0.11100541){\color[rgb]{0,0,0}\makebox(0,0)[lb]{\smash{0-way}}}%
    \put(0.08735062,0.11100541){\color[rgb]{0,0,0}\makebox(0,0)[lb]{\smash{1-way}}}%
    \put(0.1923019,0.1100766){\color[rgb]{0,0,0}\makebox(0,0)[lb]{\smash{2-way}}}%
    \put(0.03348181,0.05745283){\color[rgb]{0,0,0}\makebox(0,0)[lb]{\smash{$\mathcal{A}$}}}%
    \put(0.03255303,0.29057471){\color[rgb]{0,0,0}\makebox(0,0)[lb]{\smash{$a$}}}%
    \put(0.09849589,0.29243222){\color[rgb]{0,0,0}\makebox(0,0)[lb]{\smash{$\mathbf{a}$}}}%
    \put(0.20437596,0.29150347){\color[rgb]{0,0,0}\makebox(0,0)[lb]{\smash{$\mathbf{A}$}}}%
    \put(0.39991149,0.00375423){\color[rgb]{0,0,0}\makebox(0,0)[lb]{\smash{Tensor Algebra}}}%
    \put(0.08050091,0.33480754){\color[rgb]{0,0,0}\makebox(0,0)[lb]{\smash{Matrix Algebra}}}%
    \put(0.10964116,0.06023914){\color[rgb]{0,0,0}\makebox(0,0)[lb]{\smash{$\mathcal{A}_{k}$}}}%
    \put(0.21552123,0.0611679){\color[rgb]{0,0,0}\makebox(0,0)[lb]{\smash{$\mathcal{A}_{kl}$}}}%
    \put(0.35855219,0.06209665){\color[rgb]{0,0,0}\makebox(0,0)[lb]{\smash{$\mathcal{A}_{klm}$}}}%
    \put(0.63718395,0.0630254){\color[rgb]{0,0,0}\makebox(0,0)[lb]{\smash{$\mathcal{A}_{klmn}$}}}%
    \put(0.01119127,0.08964365){\color[rgb]{0,0,0}\makebox(0,0)[lb]{\smash{tensor}}}%
    \put(0.08735062,0.08964365){\color[rgb]{0,0,0}\makebox(0,0)[lb]{\smash{tensor}}}%
    \put(0.19137313,0.08964365){\color[rgb]{0,0,0}\makebox(0,0)[lb]{\smash{tensor}}}%
    \put(0.33626164,0.0905724){\color[rgb]{0,0,0}\makebox(0,0)[lb]{\smash{tensor}}}%
    \put(0.59538917,0.09150115){\color[rgb]{0,0,0}\makebox(0,0)[lb]{\smash{tensor}}}%
    \put(0,0){\includegraphics[width=\unitlength,page=5]{Tensor.pdf}}%
    \put(0.92446165,0.11174475){\color[rgb]{0,0,0}\makebox(0,0)[lb]{\smash{N-way}}}%
    \put(0.92944856,0.09237185){\color[rgb]{0,0,0}\makebox(0,0)[lb]{\smash{tensor}}}%
    \put(0,0){\includegraphics[width=\unitlength,page=6]{Tensor.pdf}}%
  \end{picture}%
\endgroup%

%% file: FiltAlg.pdf_tex
\begingroup%
  \makeatletter%
  \providecommand\color[2][]{%
    \errmessage{(Inkscape) Color is used for the text in Inkscape, but the package 'color.sty' is not loaded}%
    \renewcommand\color[2][]{}%
  }%
  \providecommand\transparent[1]{%
    \errmessage{(Inkscape) Transparency is used (non-zero) for the text in Inkscape, but the package 'transparent.sty' is not loaded}%
    \renewcommand\transparent[1]{}%
  }%
  \providecommand\rotatebox[2]{#2}%
  \ifx\svgwidth\undefined%
    \setlength{\unitlength}{916.04372626bp}%
    \ifx\svgscale\undefined%
      \relax%
    \else%
      \setlength{\unitlength}{\unitlength * \real{\svgscale}}%
    \fi%
  \else%
    \setlength{\unitlength}{\svgwidth}%
  \fi%
  \global\let\svgwidth\undefined%
  \global\let\svgscale\undefined%
  \makeatother%
  \begin{picture}(1,0.22187934)%
    \put(-0.00090612,0.19621197){\color[rgb]{0,0,0}\makebox(0,0)[lb]{\smash{$\mathcal{\hat{d}}_{m_1m_2m_3}$}}}%
    \put(0,0){\includegraphics[width=\unitlength,page=1]{FiltAlg.pdf}}%
    \put(0.18792087,0.13564668){\color[rgb]{0,0,0}\makebox(0,0)[lb]{\smash{$*  \hat{f}_{k_1m_1}$}}}%
    \put(0,0){\includegraphics[width=\unitlength,page=2]{FiltAlg.pdf}}%
    \put(0.38617245,0.13564668){\color[rgb]{0,0,0}\makebox(0,0)[lb]{\smash{$*  \hat{w}_{k_2m_2}$}}}%
    \put(0.58325446,0.13564668){\color[rgb]{0,0,0}\makebox(0,0)[lb]{\smash{$*  \hat{f}_{k_3m_3}$}}}%
    \put(0,0){\includegraphics[width=\unitlength,page=3]{FiltAlg.pdf}}%
    \put(0.38617245,0.07833502){\color[rgb]{0,0,0}\makebox(0,0)[lb]{\smash{$*  \hat{f}_{k_2m_2}$}}}%
    \put(0.58325446,0.07833502){\color[rgb]{0,0,0}\makebox(0,0)[lb]{\smash{$*  \hat{w}_{k_3m_3}$}}}%
    \put(0,0){\includegraphics[width=\unitlength,page=4]{FiltAlg.pdf}}%
    \put(0.18792092,0.19295834){\color[rgb]{0,0,0}\makebox(0,0)[lb]{\smash{$*  \hat{w}_{k_1m_1}$}}}%
    \put(0,0){\includegraphics[width=\unitlength,page=5]{FiltAlg.pdf}}%
    \put(0.38617245,0.19295834){\color[rgb]{0,0,0}\makebox(0,0)[lb]{\smash{$*  \hat{f}_{k_2m_2}$}}}%
    \put(0.58325446,0.19295834){\color[rgb]{0,0,0}\makebox(0,0)[lb]{\smash{$*  \hat{f}_{k_3m_3}$}}}%
    \put(0,0){\includegraphics[width=\unitlength,page=6]{FiltAlg.pdf}}%
    \put(-0.00090616,0.02427697){\color[rgb]{0,0,0}\makebox(0,0)[lb]{\smash{$\mathcal{\hat{m}}_{m_1m_2m_3}$}}}%
    \put(0,0){\includegraphics[width=\unitlength,page=7]{FiltAlg.pdf}}%
    \put(0.18792089,0.02102331){\color[rgb]{0,0,0}\makebox(0,0)[lb]{\smash{$*  \hat{f}_{k_1m_1}$}}}%
    \put(0,0){\includegraphics[width=\unitlength,page=8]{FiltAlg.pdf}}%
    \put(0.3861724,0.02102331){\color[rgb]{0,0,0}\makebox(0,0)[lb]{\smash{$*  \hat{f}_{k_2m_2}$}}}%
    \put(0.58325441,0.02102331){\color[rgb]{0,0,0}\makebox(0,0)[lb]{\smash{$*  \hat{f}_{k_3m_3}$}}}%
    \put(0,0){\includegraphics[width=\unitlength,page=9]{FiltAlg.pdf}}%
    \put(0.80712213,0.10689054){\color[rgb]{0,0,0}\makebox(0,0)[lb]{\smash{$\mathcal{\hat{p}}_{k_1k_2k_3}$}}}%
  \end{picture}%
\endgroup%

%% file: Solver.pdf_tex
\begingroup%
  \makeatletter%
  \providecommand\color[2][]{%
    \errmessage{(Inkscape) Color is used for the text in Inkscape, but the package 'color.sty' is not loaded}%
    \renewcommand\color[2][]{}%
  }%
  \providecommand\transparent[1]{%
    \errmessage{(Inkscape) Transparency is used (non-zero) for the text in Inkscape, but the package 'transparent.sty' is not loaded}%
    \renewcommand\transparent[1]{}%
  }%
  \providecommand\rotatebox[2]{#2}%
  \ifx\svgwidth\undefined%
    \setlength{\unitlength}{1030.09234318bp}%
    \ifx\svgscale\undefined%
      \relax%
    \else%
      \setlength{\unitlength}{\unitlength * \real{\svgscale}}%
    \fi%
  \else%
    \setlength{\unitlength}{\svgwidth}%
  \fi%
  \global\let\svgwidth\undefined%
  \global\let\svgscale\undefined%
  \makeatother%
  \begin{picture}(1,0.33588841)%
    \put(0,0){\includegraphics[width=\unitlength,page=1]{Solver.pdf}}%
    \put(0.04334086,0.2017894){\color[rgb]{0,0,0}\makebox(0,0)[lb]{\smash{Ritz-Galerkin\\  }}}%
    \put(0.04631446,0.18675979){\color[rgb]{0,0,0}\makebox(0,0)[lb]{\smash{formulation\\  }}}%
    \put(0.05299253,0.17077618){\color[rgb]{0,0,0}\makebox(0,0)[lb]{\smash{of a PDE\\  }}}%
    \put(0,0){\includegraphics[width=\unitlength,page=2]{Solver.pdf}}%
    \put(0.00567363,0.15214117){\color[rgb]{0,0,0}\makebox(0,0)[lb]{\smash{$a(\phi(x),\psi(x))=l(\psi(x))$\\  }}}%
    \put(0,0){\includegraphics[width=\unitlength,page=3]{Solver.pdf}}%
    \put(0.00537326,0.0980446){\color[rgb]{0,0,0}\makebox(0,0)[lb]{\smash{PDE coefficients, source\\  }}}%
    \put(0.05519161,0.08143194){\color[rgb]{0,0,0}\makebox(0,0)[lb]{\smash{$D(x),...,q(x)$\\  }}}%
    \put(0,0){\includegraphics[width=\unitlength,page=4]{Solver.pdf}}%
    \put(0.40198751,0.29161647){\color[rgb]{0,0,0}\makebox(0,0)[lb]{\smash{$\mathcal{w}_{k_1m_1},~\mathcal{w}_{k_2m_2},...$\\  }}}%
    \put(0.40468585,0.27826329){\color[rgb]{0,0,0}\makebox(0,0)[lb]{\smash{$\mathcal{f}_{k_1m_1},\mathcal{f}_{k_2m_2},...$\\  }}}%
    \put(0,0){\includegraphics[width=\unitlength,page=5]{Solver.pdf}}%
    \put(0.40373601,0.10583735){\color[rgb]{0,0,0}\makebox(0,0)[lb]{\smash{B-spline coeffs.\\  }}}%
    \put(0.40291539,0.0891797){\color[rgb]{0,0,0}\makebox(0,0)[lb]{\smash{$\mathcal{d}_{\mathbf{m}},...,\mathcal{q}_{\mathbf{l}}$\\  }}}%
    \put(0,0){\includegraphics[width=\unitlength,page=6]{Solver.pdf}}%
    \put(0.58091455,0.2144252){\color[rgb]{0,0,0}\makebox(0,0)[lb]{\smash{Solving \\  }}}%
    \put(0.58043755,0.19475454){\color[rgb]{0,0,0}\makebox(0,0)[lb]{\smash{the system  \\  }}}%
    \put(0.57948354,0.17556899){\color[rgb]{0,0,0}\makebox(0,0)[lb]{\smash{of equations  \\  }}}%
    \put(0,0){\includegraphics[width=\unitlength,page=7]{Solver.pdf}}%
    \put(0.76715084,0.18902673){\color[rgb]{0,0,0}\makebox(0,0)[lb]{\smash{Indirect \\  }}}%
    \put(0.8221971,0.18938651){\color[rgb]{0,0,0}\makebox(0,0)[lb]{\smash{B-spline \\  }}}%
    \put(0.78478,0.17499533){\color[rgb]{0,0,0}\makebox(0,0)[lb]{\smash{transform \\  }}}%
    \put(0,0){\includegraphics[width=\unitlength,page=8]{Solver.pdf}}%
    \put(0.93416133,0.19681789){\color[rgb]{0,0,0}\makebox(0,0)[lb]{\smash{$\phi(x')$\\  }}}%
    \put(0.7375882,0.11635263){\color[rgb]{0,0,0}\makebox(0,0)[lb]{\smash{$\phi(x')=\sum c_{k} \beta^{n}(x'-k)$\\  }}}%
    \put(0,0){\includegraphics[width=\unitlength,page=9]{Solver.pdf}}%
    \put(0.08500242,0.01794277){\color[rgb]{0,0,0}\makebox(0,0)[lb]{\smash{continuous domain\\  }}}%
    \put(0.24874934,0.10013511){\color[rgb]{0,0,0}\makebox(0,0)[lb]{\smash{Direct B-spline\\  }}}%
    \put(0.23040062,0.06559627){\color[rgb]{0,0,0}\makebox(0,0)[lb]{\smash{(recursive filtration)\\  }}}%
    \put(0.26170143,0.08502438){\color[rgb]{0,0,0}\makebox(0,0)[lb]{\smash{transform\\  }}}%
    \put(0.21722145,0.20547792){\color[rgb]{0,0,0}\makebox(0,0)[lb]{\smash{B-spline approximation\\  }}}%
    \put(0.2563475,0.18793867){\color[rgb]{0,0,0}\makebox(0,0)[lb]{\smash{of a solution\\  }}}%
    \put(0.21830342,0.1661726){\color[rgb]{0,0,0}\makebox(0,0)[lb]{\smash{$\phi(x)=\sum c_{k} \beta^{n}(x-k)$\\  }}}%
    \put(0,0){\includegraphics[width=\unitlength,page=10]{Solver.pdf}}%
    \put(0.85009527,0.01698877){\color[rgb]{0,0,0}\makebox(0,0)[lb]{\smash{continuous domain\\  }}}%
    \put(0.39051326,0.1941208){\color[rgb]{0,0,0}\makebox(0,0)[lb]{\smash{$\mathcal{c}_{\mathbf{k}}$\\  }}}%
    \put(0,0){\includegraphics[width=\unitlength,page=11]{Solver.pdf}}%
    \put(0.21903542,0.28326711){\color[rgb]{0,0,0}\makebox(0,0)[lb]{\smash{Compute decomposed \\  }}}%
    \put(0.21903542,0.26311945){\color[rgb]{0,0,0}\makebox(0,0)[lb]{\smash{domain and full boundary  \\  }}}%
    \put(0.21903542,0.2444109){\color[rgb]{0,0,0}\makebox(0,0)[lb]{\smash{kernels  \\  }}}%
    \put(0.69085307,0.19514698){\color[rgb]{0,0,0}\makebox(0,0)[lb]{\smash{$\mathcal{c}_{k}$\\  }}}%
    \put(0.40516285,0.26538416){\color[rgb]{0,0,0}\makebox(0,0)[lb]{\smash{$...$\\  }}}%
    \put(0,0){\includegraphics[width=\unitlength,page=12]{Solver.pdf}}%
    \put(0.53016065,0.01794277){\color[rgb]{0,0,0}\makebox(0,0)[lb]{\smash{Tensor B-spline space\\  }}}%
    \put(0.77857891,0.14589801){\color[rgb]{0,0,0}\makebox(0,0)[lb]{\smash{on grid $x'$ \\  }}}%
    \put(0.76045279,0.16068519){\color[rgb]{0,0,0}\makebox(0,0)[lb]{\smash{evaluate solution \\  }}}%
    \put(0.10834384,0.31835219){\color[rgb]{0,0,0}\makebox(0,0)[lb]{\smash{1}}}%
    \put(0,0){\includegraphics[width=\unitlength,page=13]{Solver.pdf}}%
    \put(0.29723781,0.31930618){\color[rgb]{0,0,0}\makebox(0,0)[lb]{\smash{2}}}%
    \put(0,0){\includegraphics[width=\unitlength,page=14]{Solver.pdf}}%
    \put(0.61730815,0.31930618){\color[rgb]{0,0,0}\makebox(0,0)[lb]{\smash{3}}}%
    \put(0,0){\includegraphics[width=\unitlength,page=15]{Solver.pdf}}%
    \put(0.8124032,0.32026019){\color[rgb]{0,0,0}\makebox(0,0)[lb]{\smash{4}}}%
    \put(0,0){\includegraphics[width=\unitlength,page=16]{Solver.pdf}}%
    \put(0.60182642,0.15262138){\color[rgb]{0,0,0}\makebox(0,0)[lb]{\smash{$F(c)=t$\\  }}}%
    \put(0.54230155,0.2909419){\color[rgb]{0,0,0}\makebox(0,0)[lb]{\smash{\\  }}}%
  \end{picture}%
\endgroup%